\documentclass[12pt]{article}
\usepackage{amsmath,amssymb,hhline}
\usepackage[russian]{babel}

\usepackage[T2A]{fontenc}    
\usepackage[cp1251]{inputenc}
\usepackage[russian]{babel} 

\begin{document}

\begin{center}
{\large\bf Р.~М.~Тригуб \qquad\qquad R.~M.~Trigub}
\end{center}

\begin{center}
{\large \bf О разных модулях гладкости и $K$-функционалах\\ On various moduli of smoothness and $K$-functionals}
\end{center}

\emph{\textbf{Аннотация.}} Статья посвящена определению точного порядка приближения
функций линейными средними рядов и интегралов Фурье и определению соответствующих
$K$-функционалов через специальные модули гладкости.

$\S 1$ --- введение, а в $\S2$ изучаются функции на прямой $\mathbb{R}$. Типичный
результат, полученный давно: для любой $2\pi$-периодической функции из $L_p$ на периоде
при любом $p\in[1,+\infty]$ $(L_\infty=C)$, при любом $r\in\mathbb{N}$ существует
тригонометрический полином $\tau_{r,n}(f)$ порядка не выше $n$ такой, что
\begin{equation*}
    \big\|f-\tau_{r,n}(f)\big\|_p\asymp\omega_r\Big(f;\frac{1}{n}\Big)_p\asymp\inf\limits_{g}\Big\{\|f-g\|_p+\frac{1}{n^r}\big\|g^{(r)}\big\|_p\Big\}
\end{equation*}
(двусторонние или порядковые неравенства с положительными константами, зависящими лишь
от $r$). В $\S 3$ --- функции на $\mathbb{R}^d$ ($d\geq2$), а в $\S 4$ --- в банаховых
пространствах.

Статья носит частично обзорный характер. Доказательства приведены лишь для теорем 2.2,
3.9 и в $\S 4$. В $\S 5$ сформулированы некоторые примыкающие нерешенные вопросы.

Список литературы: 52 названия.

\emph{\textbf{Abstract.}} In this paper, exact rate of approximation of functions by
linear means of Fourier series and Fourier integrals and corresponding $K$-functionals
are expressed via special moduli of smoothness.  .

Introduction is given in $\S 1$. In $\S2$ functions on the line $\mathbb{R}$ are studied.
A typical (well-known) result is as follows: for each $2\pi$-periodic function in $L_p$ 
on the period, for any $p\in[1,+\infty]$ ($L_\infty=C$) and $r\in\mathbb{N}$, there is
a trigonometric polynomial $\tau_{r,n}(f)$ of degree not greater than $n$ such that
\begin{equation*}
    \big\|f-\tau_{r,n}(f)\big\|_p\asymp\omega_r\Big(f;\frac{1}{n}\Big)_p\asymp
    \inf\limits_{g}\Big\{\|f-g\|_p+\frac{1}{n^r}\big\|g^{(r)}\big\|_p\Big\},
\end{equation*}
where the positive constants in these bilateral inequalities depend only on $r$. In $\S 3$
we deal with functions on $\mathbb{R}^d$ ($d\geq2$), while in $\S 4$ with functions on
Banach spaces. 

The paper is partially of survey nature. The proofs are given only for Theorems 2.2,
3.9 and those in $\S 4$. Related open problems are formulated in $\S 5$. 

The list of references contains 52 items.

\emph{\textbf{Ключевые слова:}} ряд и интеграл Фурье, модуль гладкости и $K$-функционал,
винеровские банаховы алгебры, принцип сравнения мультипликаторов Фурье.



\newpage

\begin{center}
    \textbf{\S1 Введение}
\end{center}

Модуль гладкости непрерывной $2\pi$-периодической функции
$f:\mathbb{R}\rightarrow\mathbb{C}\ $ ($f\in C(\mathbb{T})$, $\mathbb{T}=[-\pi,\pi]$)
порядка $r\in\mathbb{N}$ и шага $h>0$ определяют так:
\begin{equation*}
    \omega_r(f;h)=\sup\Big\{\big|\Delta^r_\delta f(x)\big|:
    x\in\mathbb{T},\delta\in(0,h]\Big\},
\end{equation*}
\begin{equation*}
\Delta^r_\delta f(x)=\sum\limits_{\nu=0}^r\binom{r}{\nu}(-1)^\nu f(x+\nu\delta).
\end{equation*}

Модуль непрерывности $\omega=\omega_1$ использовал Лебег. При $r\geq2$ модуль гладкости
ввел С.~Н.~Бернштейн (1923), а основные его свойства изучил Маршо (A. Marchaud, 1927).

Известно, напр., что
\begin{equation*}
    \omega_r(f;h)=O(h^\alpha)\qquad (0<\alpha\leq r,\ h\rightarrow0)
\end{equation*}
при нецелом $\alpha$ тогда и только тогда, когда $f^{[\alpha]}\in {\rm Lip}\
\big(\alpha-[\alpha]\big)$, при целом $\alpha<r$, если
$\omega_2\big(f^{(\alpha-1)};h\big)=O(h)$ и при $\alpha=r$ --- $f^{(r-1)}\in {\rm Lip}\
1$.

Подобным образом определяется модуль гладкости функции $\mathbb{R}\rightarrow\mathbb{C}$
или на отрезке прямой $\mathbb{R}$. И для функций из $L_p$.

Теорему о продолжении функции с сохранением порядка убывания $\omega_r$ при
$h\rightarrow0$ получил О.~В.~Бесов [\ref{Besov}]. Некоторые свойства (о модуле
гладкости произведения функций, напр.) получил автор (см. [\ref{Trigub_Belinsky}]). См.
также [\ref{De_Vore}], [\ref{Dzyadyk_Shevchuk}].

В случае числовых функций нескольких переменных имеются разные модули гладкости (полные,
частные, смешанные). По поводу соотношений между ними см. М.~Ф.~Тиман [\ref{Timan}], а
общую теорему о продолжении --- см. Ю.~А.~Брудный [\ref{Brudnyi}].

С.~Н.~Бернштейн [\ref{Bernstein}] рассмотрел модуль гладкости, определяемый следующей
"разностью"\ при $r\geq2$
\begin{equation*}
\begin{split}
    f(x)-\sum\limits_{k=1}^ra_k f(x+b_k\delta),\quad &\sum\limits_{k=1}^ra_k=1,\quad
    \sum\limits_{k=1}^ra_kb_k^s=0\quad (1\leq s\leq r-1),\\
    &1=b_1<b_2<...<b_r.
\end{split}
\end{equation*}

Это частный случай такого модуля ($f$ --- ограничена и непрерывна на прямой).
\begin{equation*}
    \sup\limits_x\sup\limits_{0<\delta\leq h}\Bigg|\int\limits_{-\infty}^\infty
    f(x-\delta u)d\mu(u)\Bigg|,\qquad \int\limits_{-\infty}^\infty
    u^kd\mu(u)du=0\quad (0\leq k\leq r-1),
\end{equation*}
который изучали H.~Shapiro и J.~Boman (см., напр., [\ref{Boman}]).

Ранее [\ref{Trigub1968}] автор ввел следующий модуль
\begin{equation*}
    \tilde{\omega}_r(f;h)=\sup\limits_x\Bigg|\frac{1}{h}\int\limits_0^h\Delta_\delta^r
    f(x)d\delta\Bigg|
\end{equation*}
(верхняя грань по шагу заменена интегральным средним), который по порядку эквивалентен
$\omega_r(f;h)$. Этот линеаризованный модуль удобен в применении
([\ref{Trigub_Belinsky}], см. ниже в $\S2$).

В $60^\text{е}$-годы прошлого века для построения интерполяционных пространств между
двумя банаховыми пространствами $E_1$ и $E_2$ вещественным методом были введены
(J.~L.~Lions и J.~Peetre) $K$-функционалы ($f\in E_1+E_2$, $h>0$)
\begin{equation*}
    K(f;h)=K(f;h,E_1,E_2)=\inf\Big\{\|f_1\|_{E_1}+h\|f_2\|_{E_2},\ f=f_1+f_2\Big\}.
\end{equation*}

См. [\ref{Bergh}], а также [\ref{De_Vore}].

В случае, когда одно из пространств является пространством гладких функций,
$K$-функционал, как и модуль гладкости, определяет при $h\rightarrow0$ промежуточную
гладкость.

Z.~Ciesielski (1983) [\ref{Ciesielski}] поставил вопрос о формулах для $K$-функционалов,
выражаемых для периодических функций через средние их рядов Фурье.

В работе [\ref{Trigub1988}] найдены такие формулы для пространств, определяемых
дифференциальными операторами $\Delta^r$ ($\Delta$ --- оператор Лапласа) и
$D_{2r}=\sum\limits_{j=1}^d\frac{\partial^{2r}}{\partial x_j^{2r}}$ ($d\in\mathbb{N}$,
$x=(x_1,...,x_d)$, $r\in\mathbb{N}$). В случае оператора Лапласа $\Delta^1$ другая форму
указана в [\ref{Ditzian}]. См. ниже 3.4--3.5. При этом автор использовал двусторонние
(порядковые) неравенства для приближения индивидуальных функций классическими средними
рядов Фурье (одномерный случай см., напр., в [\ref{Trigub1965}] (1965)). Этой проблемой
успешно занимались В.~В.~Жук ([\ref{Zhuk1967}], см. также [\ref{Zhuk1982}]),
Э.~В.~Стороженко [\ref{Storozhenko}], М.~Ф.~Тиман и В.~П.~Пономаренко
[\ref{Timan_Ponomarenko}].

А в случае любого числа переменных --- Э.~С.~Белинский [\ref{Belinsky}], автор
[\ref{Trigub1980}], О.~И.~Кузнецова [\ref{Kuznetsova_Trigub}], Ю.~Л.~Носенко
[\ref{Nosenko}]. А после статьи Z.~Ditzian--K.~G.~Ivanov [\ref{Ditzian_Ivanov}] (1993)
такие результаты о точном порядке приближения функций известными полиномами часто
называют "strong converse inequalities". Назову еще V.~Totik, X.~L.~Zhou,
К.~В.~Руновский и Ю.~С.~Коломойцев. См. ссылки в статье B.~R.~Draganov [\ref{Draganov}].
Еще нужно упомянуть о подобных результатах в пространстве $L_p$, $p>0$, и пространствах
Харди $H_p$, $0<p\leq1$, на диске [\ref{Trigub1997}], полуплоскости [\ref{Tovstolis}] и
шаре [\ref{Vit.Volchkov}].

Для определения точного порядка приближения функции на торе $\mathbb{T}^d$ или на всем
пространстве $\mathbb{R}^d$ используют обычно один из следующих двух методов. Первый из
них основан на прямых теоремах теории приближений (типа Джексона) и экстремальных
свойствах целых функций экспоненциального типа (типа неравенства Бернштейна). Этот метод
применим и к нелинейным операторам приближения. Второй основан на принципе сравнения
мультипликаторов (впервые применен [\ref{Shapiro}] и [\ref{Trigub1968}]). Таким образом
получаются и прямые теоремы, и экстремальные свойства полиномов и целых функций
экспоненциального типа. Достоинство этого метода: точность.

О мультипликаторах Фурье см. [\ref{Stein}] ($L_p$, $p\in(1,+\infty)$) и
[\ref{Stein_Weiss}] ($p=1$ и $p=\infty$). О связи мультипликаторов Фурье для интегралов
и рядов Фурье см. [\ref{Stein_Weiss}] (гл.~VII).

Мультипликатор Фурье интегралов Фурье из $L_1$ в $L_1$ определяется преобразованием
Фурье конечной комплекснозначной борелевской меры (по поводу свойств таких мер см.
[\ref{Makarov_Podkorytov}]).
\begin{equation*}
    B(\mathbb{R}^d)=\Big\{f(x)=\int\limits_{\mathbb{R}^d}e^{-i(x,u)}d\mu(u),\
    \|f\|_B={\rm var}\mu<\infty\Big\}
\end{equation*}
\begin{equation*}
    A(\mathbb{R}^d)=\Big\{f(x)=\hat{g}(x)=\int\limits_{\mathbb{R}^d}e^{-i(x,u)}g(u)du,\
    \|f\|_A=\|g\|_{L_1}<\infty\Big\}
\end{equation*}

Этим банаховым алгебрам посвящена обзорная статья [\ref{Liflyand_Samko_Trigub}].

Для перехода от приближений к $K$-функционалам применяем простую лемму (см. ниже в конце
введения).

Таким образом, нелинейная задача определения $K$-функционала решается применением
линейных операторов (и даже сверточных). См. также [\ref{Trigub1991}--\ref{Trigub2005}].

Ряд Фурье $2\pi$-периодической по $x_j$ ($1\leq j\leq d$) функции $f\in
L_1(\mathbb{T}^d)$, $\mathbb{T}^d=[-\pi,\pi]^d$ будем писать в виде
$\big((x,y)=\sum\limits_{j=1}^d x_jy_j$, $|x|=\sqrt{(x,x)}\big)$
\begin{equation*}
    f\sim \sum\limits_{k\subset \mathbb{Z}^d}\hat{f}_ke_k,\quad e_k=e^{i(k,x)},\quad
    \hat{f}_k=\frac{1}{(2\pi)^d}\int\limits_{\mathbb{T}^d}f(x)e^{-i(k,x)}dx.
\end{equation*}

Средние $\Phi_\varepsilon(f)$ рядов Фурье, определяемые функцией
$\varphi:\mathbb{R}^d\rightarrow\mathbb{C}$, имеют вид ($\varepsilon>0$)
\begin{equation*}
\Phi_\varepsilon(f)=\Phi_\varepsilon(f,x)\sim\sum\limits_{k\in\mathbb{Z}^d}\varphi(\varepsilon
k)\hat{f}_ke_k.
\end{equation*}

\textbf{Принцип сравнения.}

\emph{Если $\varphi$ и $\psi:\mathbb{R}^d\rightarrow\mathbb{C}$ ограничены и непрерывны
почти всюду, а переходная функция $g=\frac{1-\varphi}{1-\psi}$ после исправления
(доопределения) по непрерывности принадлежит $B(\mathbb{R}^d)$, то при любом $p\geq1$ и
$\varepsilon>0$
\begin{equation*}
    \Big\|f-\sum\limits_{k\in\mathbb{Z}^d}\varphi(\varepsilon
    k)\hat{f}_ke_k\Big\|_p\leq\|g\|_B\cdot\Big\|f-\sum\limits_{k\in\mathbb{Z}^d}\psi(\varepsilon
    k)\hat{f}_ke_k\Big\|_p
\end{equation*}
для всех функций $f\in L_p(\mathbb{T}^d)$, для которых конечна правая часть} (см.
[\ref{Trigub_Belinsky}], \textbf{7.1.11}, \textbf{7.1.12} и замечание перед следствием
1).

Напомним определение $K$-функционала.

Пусть $d_\alpha$ --- дифференциальный оператор, определяемый матрицей
$\{\mu_{k,\alpha}\}_{k\in\mathbb{Z}^d\setminus\{0\}}$ такой, что $\mu_{k,\alpha}\neq0$ и
$\lim\limits_{|k|\rightarrow\infty}\mu_{k,\alpha}=\infty$. А именно
\begin{equation*}
    d_\alpha f\sim\sum\limits_{k\neq0}\mu_{k,\alpha}\hat{f}_ke_k
\end{equation*}
(при $\mu_{k,\alpha}=-|k|^{2\alpha}$ это $\Delta^\alpha$, где $\Delta$ --- оператор
Лапласа).

Если
\begin{equation*}
    W(d_\alpha)_p=\Big\{f\in L_p(\mathbb{T}^d): d_\alpha f\in L_p(\mathbb{T}^d)\Big\},
\end{equation*}
то
\begin{equation*}
    K\big(\varepsilon;f,L_p,W(d_\alpha)_p\big)=\inf\limits_g\Big\{\|f-g\|_p+\varepsilon\|d_\alpha
    g\|_p\Big\}.
\end{equation*}

Для перехода к $K$-функционалам используется следующая лемма (доказательство приведено в
$\S 3$).

\textbf{Лемма.}

\emph{Пусть $d_\alpha \Phi_\varepsilon(f)\in L_p$, $p\geq1$, $\varepsilon>0$. Если при
некотором $a$, не зависящем от $f$ и $\varepsilon$, выполняются следующие три условия
($\alpha)$ и $\beta)$ и необходимы)
\begin{equation*}
\begin{split}
    &\alpha)\ \big\|f-\Phi_\varepsilon(f)\big\|_p\leq a\varepsilon\big\|d_\alpha
    f\big\|_p \quad (\text{в предположении, что\ } d_\alpha f\in L_p),\\
&\beta)\ \sup\limits_{\varepsilon>0}\big\|\Phi_\varepsilon (f)\big\|_p\leq a\|f\|_p,\\
&\gamma)\ \varepsilon\big\|d_\alpha\Phi_\varepsilon(f)\big\|_p\leq a
\big\|f-\Phi_\varepsilon(f)\big\|_p,
\end{split}
\end{equation*}
то
\begin{equation*}
K\big(\varepsilon;f,L_p,W(d_\alpha)_p\big)\asymp \big\|f-\Phi_\varepsilon(f)\big\|_p
\end{equation*}
 (двустороннее неравенство с положительными константами, не зависящими от $f$ и $\varepsilon$).}

В $\S2$ настоящей статьи изучаются функции на $\mathbb{R}$.

В $\S 3$ --- на $\mathbb{R}^d$ ($d\geq2$), а в $\S 4$ --- модули гладкости в банаховых
пространствах.

В $\S 5$ приведены примыкающие нерешенные вопросы.

Через $c(a,b)$ будем обозначать положительные величины, зависящие лишь от $a$ и $b$.

\begin{center}
    \textbf{\S2 Модули гладкости и $K$-функционалы на $\mathbb{R}$}
\end{center}

Для того чтобы объединить периодические и непериодические функции на $\mathbb{R}^1$
можно ввести семейство норм ($\mathbb{N}>0$)
\begin{equation*}
    \|f\|_{p,N}=\sup\limits_{x\in\mathbb{R}}\Bigg(\int\limits_{x-N}^{x+N}|f(u)|^pdu\Bigg)^{\frac{1}{p}}.
\end{equation*}

Если $f\in L_p(\mathbb{T})$ и $2\pi$-периодична, то $\|f\|_p=\|f\|_{p,\pi}$. А если
$f\in L_p(\mathbb{R})$, то $\|f\|_p=\lim\limits_{N\rightarrow\infty}\|f\|_{p,N}$.

Обозначим через $\tau_t$ оператор сдвига $\big((\tau_t f)(x)=f(x+t)\big)$.

\textbf{2.1.} \emph{Пусть $E$ --- линейное множество ограниченных и равномерно
непрерывных функций на $\mathbb{R}$, замкнутое относительно равномерной сходимости, а
линейный оператор $A$ действует из $E$ в $E$ и $\|Af\|_\infty\leq a\|f\|_\infty$
($\sup$-норма). Если еще множество $E$ инвариантно относительно сдвига и оператор $A$
коммутирует со сдвигом ($A\tau_t=\tau_t A$, для всех $t\in \mathbb{R}$), то и для любого
$p\in[1,+\infty)$, $N>0$ и $f\in E$
\begin{equation*}
    \|Af\|_{p,N}\leq a\|f\|_{p,N}.
\end{equation*}
} Доказательство см. [\ref{Trigub_Belinsky}], \textbf{1.2.7}.

Введем теперь линеаризованный модуль гладкости ($f\in L_{1,loc}(\mathbb{R})$):
\begin{equation*}
    \tilde{\omega}_r(f;h)=\Bigg\|\frac{1}{h}\int\limits_0^h\Delta_\delta^r
    f(\cdot)d\delta\Bigg\|=\Bigg\|\int\limits_0^1\Delta_{th}^r
    f(\cdot)dt\Bigg\|
\end{equation*}
(верхняя грань по шагу $\delta$ заменена интегральным средним).

Очевидно, что при любом $h>0\ $ $\tilde{\omega}_r(f;h)\leq \omega_r(f;h)$.

\textbf{2.2.} \emph{Для любого $p\in[1,+\infty]$ и любой функции $f\in L_p$
($L_p(\mathbb{R})$ или $L_p(\mathbb{T})$, если функция периодическая) при всех
достаточно малых $h>0$ (в зависимости от $r\in\mathbb{N}$)
\begin{equation*}
    \omega_r(f;h)_p\leq c(r)\tilde{\omega}_r(f;h)_p.
\end{equation*}
}

$\rhd$ Нужно доказать, что
\begin{equation*}
    \sup\limits_{\theta\in (0,1]}\big\|\Delta_{\theta h}^r f(\cdot)\big\|\leq
    c(r)\tilde{\omega}_r(f;h).
\end{equation*}

Применим принцип сравнения (см. выше), учитывая, что
\begin{equation*}
\Delta_{\theta h}^r f\sim\sum\limits_{k\in \mathbb{Z}}\big(1-e^{i\theta
hk}\big)^r\hat{f}_ke_k,
\end{equation*}
\begin{equation*}
    \int\limits_0^1\Delta_{th}^r f(\cdot)dt=\sum\limits_{k\in
    \mathbb{Z}}\int\limits_0^1\big(1-e^{ithk}\big)^rdt\hat{f}_ke_k.
\end{equation*}

Нужно доказать, что
\begin{equation*}
    \begin{split}
        & \sup\limits_\theta\|g_{r,\theta}\|_B<\infty,\quad g_{r,\theta}=\frac{\varphi_{r,\theta}}{\psi_r},\quad \varphi_{r,\theta}(x)=\big(1-e^{i\theta x}\big)^r, \\
        & \psi_r(x)=\int\limits_0^1\big(1-e^{it
        x}\big)^rdt=1+\sum\limits_{\nu=1}^r(-1)^\nu\binom{r}{\nu}\frac{e^{i\nu x}-1}{i\nu
        x}.
    \end{split}
\end{equation*}

Имеем $g_{r,\theta}(0)=(r+1)\theta^r$. А то, что $\psi_r(x)\neq0$ при $x\in
\mathbb{R}\setminus\{0\}$ доказано автором с использованием теоремы Линдемана о
трансцендентности значений показательной функции (см. лемму \textbf{8.3.5 в)} в
[\ref{Trigub_Belinsky}]).

Учитывая, что при $\lambda\in \mathbb{R}\ $ $\big\|e^{i\lambda(\cdot)}\big\|_B=1$,
получаем
\begin{equation*}
    \|g_{r,\theta}\|_B=\Big\|\frac{\varphi_{r,\theta}(1-\psi_r)}{\psi_r}+\varphi_r\Big\|_B\leq\Big\|\frac{\varphi_{r,\theta}(1-\psi_r)}{\psi_r}\Big\|_A+2^r.
\end{equation*}

Приведем некоторые достаточные условия принадлежности $A(\mathbb{R})$.

\textbf{А$_1$.} \textbf{Аналог признака Бернштейна-Саса абсолютной сходимости рядов
Фурье.}

\emph{Если при $p\in (1,2]$ $f\in L_p(\mathbb{R})\cap C_0(\mathbb{R})$ и при
$\alpha\in\Big(\frac{1}{p},1\Big]$
\begin{equation*}
    \omega(f;h)_p\leq h^\alpha\qquad (h\in(0,1]),
\end{equation*}
то $f\in A(\mathbb{R})$} (E.~C.~Titchmarsh -- 1927). См. \textbf{5.1} в
[\ref{Liflyand_Samko_Trigub}].

\textbf{Б$_1$.} \emph{Если $f\in C_0(\mathbb{R})\cap L_{1,loc}(\mathbb{R})$, а
\begin{equation*}
    \int\limits_{0}^\infty\underset{|u|\geq x}{{\rm ess}\sup}\big|f'(u)\big|dx<\infty,
\end{equation*}
то $f\in A(\mathbb{R})$ тогда и только тогда, когда
\begin{equation*}
    \int\limits_{0}^\infty\frac{\big|f(x)-f(-x)\big|}{x}dx<\infty
\end{equation*}}
(см. [\ref{Trigub_Belinsky}], \textbf{6.5.9} и следствие \textbf{12.7} в
[\ref{Liflyand_Samko_Trigub}]).

\textbf{В$_1$.} \emph{Если $f\in C^1(\mathbb{R})$, а при $|x|\rightarrow\infty$
\begin{equation*}
    f(x)=O\Big(\frac{1}{|x|^\alpha}\Big)\quad (\alpha>0), \qquad
    f'(x)=O\Big(\frac{1}{|x|^\beta}\Big)\quad (\beta\in\mathbb{R}),
\end{equation*}
то при $\alpha+\beta>1$ $f\in A(\mathbb{R})$ } ([\ref{Liflyand_Trigub}], см. также
следствие \textbf{10.6} в [\ref{Liflyand_Samko_Trigub}]).

Очевидно, что равномерно по $\theta\in(0,1]$ при $|x|\rightarrow\infty$
\begin{equation*}
\frac{\varphi_{r,\theta}(x)\big(1-\psi_r(x)\big)}{\psi_r(x)}=O\Big(\frac{1}{|x|}\Big).
\end{equation*}

А так как
\begin{equation*}
    \sup\limits_x\Big|\frac{\varphi_{r,\theta}(x)}{x^r}\Big|=\sup\limits_x\Big|\frac{1-e^{i\theta
    x}}{x}\Big|^r\leq\theta^r\leq1,\qquad
    \sup\limits_x\Big|\frac{\psi_{r}(x)}{x^r}\Big|\leq\frac{1}{r+1}
\end{equation*}
и по неравенству Бернштейна для производных целых функций экспоненциального типа
\begin{equation*}
    \sup\limits_x\Big|\Big(\frac{\varphi_{r,\theta}(x)}{x^r}\Big)'\Big|\leq r\theta^r\leq r,\qquad
    \sup\limits_x\Big|\Big(\frac{\psi_{r}(x)}{x^r}\Big)'\Big|\leq\frac{r}{r+1},
\end{equation*}
то и производная
\begin{equation*}
    \Bigg(\frac{\varphi_{r,\theta}(x)\big(1-\psi_r(x)\big)}{\psi_r(x)}\Bigg)'=\Bigg(\frac{\frac{\varphi_{r,\theta}(x)}{x^r}\big(1-\psi_r(x)\big)}{\frac{\psi_r(x)}{x^r}}\Bigg)'=O\Big(\frac{1}{|x|}\Big)
\end{equation*}
равномерно по $\theta\in(0,1]$.

И можно применить $A_1$ при $p=2$ $\big(\omega(f;h)_2\leq\|f'\|_2h\big)$ или $B_1$
$(\alpha=\beta=1)$.\qquad\qquad$\blacktriangleleft$

\textbf{2.3.} \emph{Если при $\varepsilon>0$ и $n=\Big[\frac{1}{\varepsilon}\Big]$
\begin{equation*}
    \tau_{r,n}(f)=\sum\limits_{k\in\mathbb{Z}}\varphi_r(k\varepsilon)\hat{f}_ke_k,
\end{equation*}
где $\varphi_r(x)=\big(1-|x|^r\big)_+$ для четного $r$ и
$\varphi_r(x)=\big(1-|x|^{r+1}\big)_++i|x|^r\big(1-|x|\big)_+{\rm sign}~x$ для нечетного
$r$, то при $p\in[1,+\infty]$
\begin{equation*}
    \big\|f-\tau_{r,n}(f)\big\|_p\asymp\omega_r\Big(f;\frac{1}{n}\Big)_p
\end{equation*}
(двустороннее неравенство с константами, зависящими лишь от $r$).}

Доказательство аналогично предыдущему (после замены $\omega_r$ на $\tilde{\omega}_r$). А
после применения леммы, приведенной выше, получаем давно известный факт (см., напр.,
[\ref{Bergh}]):

если $d_r(f;x)=\frac{d^r}{dx^r}f(x)$, то
\begin{equation*}
    K\big(\varepsilon,f,L_p,W(d_r)_p\big)\asymp\omega_r(f;\varepsilon)_p.
\end{equation*}

\textbf{2.4.} \emph{При натуральных $\alpha$ и $r>\frac{\alpha}{2}$, $\dot{\Delta}_h
f(x)=f(x-h)-f(x+h)$
\begin{equation*}
\begin{split}
    &\Bigg\|\int\limits_{|u|\geq1}\frac{1}{|u|^{r+\alpha}}\dot{\Delta}_{\varepsilon
    u}^{2r}f(\cdot)du\Bigg\|_p\asymp\omega_\alpha(f;\varepsilon)_p\qquad (\alpha
    \text{--- четное}),\\
&\Bigg\|\int\limits_{|u|\geq1}\frac{1}{|u|^{r+\alpha}}\dot{\Delta}_{\varepsilon
    u}^{2r}f(\cdot)du\Bigg\|_p\asymp\omega_{\alpha+1}(f;\varepsilon)_p+\frac{1}{\varepsilon}\omega_{\alpha+1}(\tilde{F};\varepsilon)_p\qquad (\alpha
    \text{--- нечетное}).
    \end{split}
\end{equation*}
}

Здесь норма в $L_p(\mathbb{T})$, а $\tilde{F}$ --- сопряженная функция к первообразной
$f$.

$\rhd$ Как следует из теоремы 3.9, доказанной ниже, при $d=1$ левая часть совпадает по
порядку с
\begin{equation*}
    \Big\|f-\sum\limits_{k\in\mathbb{Z}}\big(1-(\varepsilon|k|)^\alpha\big)_+\hat{f}_ke_k\Big\|,
\end{equation*}
а точный порядок приближения такими полиномами давно известен ($\alpha$ --- четное ---
автор, $\alpha$ --- нечетное --- В.~В.~Жук, см. также [\ref{Trigub_Belinsky}],
\textbf{8.5.8}).

Вид $K$-функционала приведен ниже.

Поскольку с ростом $\alpha$ модуль $\omega_\alpha$ может только уменьшаться, то и
интегралы в левой части с ростом $\alpha$ почти убывают, если только $\alpha<2r$,
конечно.\qquad\qquad\qquad\qquad\qquad$\blacktriangleleft$

Введем теперь линеаризованный модуль гладкости нецелого порядка $r>0$.

Для $f\in C(\mathbb{T})$ и целого $q>\frac{r}{2}>0$
\begin{equation*}
    \tilde{\omega}_r(f;h)_\infty=\Bigg\|\int\limits_1^\infty\frac{1}{u^{r+1}}\Big(\dot{\Delta}_{hu}^{2q}f(\cdot)+\gamma_0\dot{\Delta}_{hu}^{2q+1}f(\cdot)\Big)du\Bigg\|_\infty
\end{equation*}
при
\begin{equation*}
    \gamma_0=\frac{1}{2}{\rm
    tg}\frac{r\pi}{2}\Bigg(\int\limits_0^\infty\frac{\sin^{2q+1}tdt}{t^{r+1}}\Bigg)^{-1}\int\limits_0^\infty\frac{\sin^{2q}t}{t^{r+1}}dt.
\end{equation*}

\textbf{2.5.} \emph{Пусть
$d_r(f)=f^{(r)}\sim\sum\limits_{k\neq0}|k|^re^{i\frac{\pi}{2}r\cdot{\rm
sign}k}\hat{f}_ke_k$, $\ r\notin\mathbb{N}$,
\begin{equation*}
\varphi_r(x)=\big(1-|x|^r\big)_+-i{\rm
tg}\Big(\frac{r\pi}{2}\Big)|x|^r\big(1-|x|\big)_+{\rm sign}~x.
\end{equation*}}
\emph{Тогда
\begin{equation*}
    K\big(\varepsilon^r,f;C,f^{(r)}\big)\asymp\Big\|f-\sum\limits_k\varphi_r(\varepsilon
    k)\hat{f}_ke_k\Big\|_\infty\asymp\tilde{\omega}_r(f;\varepsilon)_\infty
\end{equation*}
(двойные неравенства с положительными константами, зависящими лишь от $r$ и $q$).}

Эти соотношения справедливы, конечно, и в $L_p(\mathbb{T})$ при $p\in[1,+\infty)$. Но
при $p\in(1,+\infty)$ $\gamma_0$ можно выбирать произвольно (напр., $\gamma_0=0$). Тогда
константы зависят и от $p$.

См. [\ref{Trigub_Belinsky}], \textbf{8.3.1}. Здесь добавлен только $K$-функционал.

Для всех классических методов суммирования найдены и уже довольно давно точные порядки
приближения через обычные модули гладкости или специальные. Оставался только следующий
метод, введенный еще С.~Н.~Бернштейном [\ref{Bernstein}] ($D_n$ --- ядро Дирихле,
$\varphi_n(0)=1$, $s\geq r+2$)
\begin{equation*}
    \tau_{s,r,n}(f)=\gamma\int\limits_{-\pi}^\pi\sum\limits_{\nu=1}^r(-1)^{\nu+1}\binom{r}{\nu}f(x-u)D_n^s(u)du=\sum\limits_{k\in\mathbb{Z}}\varphi_n\Big(\frac{k}{n}\Big)\hat{f}_ke_k
\end{equation*}
(общая оценка приближения сверху через $\omega_r$ при $s=r+2$ доказана С.~Б.~Стечкиным
[\ref{Kolomoitsev_Trigub}], см. также [\ref{De_Vore}], [\ref{Trigub_Belinsky}]).

На этот вопрос обратил внимание В.~И.~Иванов (ТулГУ).

Как обнаружила аспирантка О.~В.~Котова (ДонНУ), уже при $r=2$ и $s=4$ $\varphi_n(x)=1$
не только при $x=0$, как это у всех классических методов суммирования рядов Фурье.

Для решения этой задачи автор ввел специальный модуль гладкости (и даже два: см.
[\ref{Kolomoitsev_Trigub}], [\ref{Trigub2013}]), а Ю.~С.~Коломойцев --- вид
$K$-функционала [\ref{Kolomoitsev_Trigub}].

Для точки $x_0=0$ полагаем
\begin{equation*}
    \tilde{\Delta}_{h,0}^rf(x)=\int\limits_0^1\Delta_{th}^rf(x)dt,
\end{equation*}
а для $x_1\in(-\pi,0)\cup(0,\pi)$
\begin{equation*}
    \tilde{\tilde{\Delta}}_{h,x_1}^1f=\int\limits_0^1\Big[\Delta_{th}^1f(x)-\lambda_1\Delta_{th}^2f(x)\Big]dt,
\end{equation*}
где
\begin{equation*}
    \lambda_1=\frac{2(ix_1+1-e^{ix_1})}{2ix_1+3-4e^{ix_1}+e^{2ix_1}}
\end{equation*}
(вещественная часть знаменателя $2(1-\cos x_1)^2\neq0$).

Пусть $X=\{x_j\}_0^m$ ($x_0=0$, $m\geq1$) --- множество различных точек из $(-\pi,\pi)$,
а $R=\{r_j\}_0^m$ --- некий набор натуральных чисел. Положим
\begin{equation*}
    K_{X,R}(\varepsilon,f)=\inf\limits_{g\in
    W_\infty^r}\Big\{\|f-g\|+\Big\|\prod\limits_{j=0}^m\Big(\varepsilon\frac{d}{dx}-ix_j\Big)^{r_j}g\Big\|\Big\},
\end{equation*}
где $r=\sum\limits_{j=0}^mr_j$.

\textbf{2.6.} \emph{Пусть $X$ и $R$ --- определенные выше наборы чисел и $m\geq1$, а
$\varphi\in C(\mathbb{R})$, ${\rm supp}\varphi\subset[-\pi,\pi]$, $\varphi(x_j)=1$
$(0\leq j\leq m)$ и только в этих точках. Если еще $\varphi^{(r_j)}(x_j)\neq0$ $(0\leq
j\leq m)$, то для того чтобы
\begin{equation*}
    K_{X,R}(\varepsilon,f)\asymp\Big\|\tilde{\Delta}_{\varepsilon,x_0}\tilde{\tilde{\Delta}}_{\varepsilon,x_1}^{r_1}\cdot\cdot\cdot\tilde{\tilde{\Delta}}_{\varepsilon,x_m}^{r_m}f\Big\|
\end{equation*}
(справа --- произведение разностных операторов, связанных с каждой точкой данного
множества) с константами, не зависящими от $f$ и $\varepsilon$, необходимо и
достаточно:
\begin{equation*}
    \frac{1-\varphi(x)}{x^r\prod\limits_{j=1}^m(x-x_j)^{r_j}}\in A(\mathbb{R}).
\end{equation*}
}

В доказательстве [\ref{Kolomoitsev_Trigub}] используются и необходимые условия в
принципе сравнения из [\ref{Trigub_Belinsky}], \textbf{7.1.11}.

\textbf{2.7.} [\ref{Trigub2013}] \emph{Пусть $s\geq6$, $2\leq r\leq s-2$ и
$r_1=2\Big[\frac{r+1}{2}\Big]$. Существует число $c(r)$  такое, что при $n\geq c(s)$ и
$h_n=\frac{2\pi}{(2n+1)s}$ для любой $f\in C(\mathbb{T})$ (определение $\tau_{s,r,n}$
см. выше)
\begin{equation*}
    \big\|f-\tau_{s,r,n}(f)\big\|\asymp\Bigg\|\tilde{\Delta}_{h_n,0}^{r_1-2p}\prod\limits_{j=1}^p\tilde{\tilde{\Delta}}_{h_n,x_{j,n}}\tilde{\tilde{\Delta}}_{h_n,-x_{j,n}}f(\cdot)\Bigg\|
\end{equation*}
(двойное неравенство с положительными константами, не зависящими от $f$ и $n$), где
$\{x_{j,n}\}_1^p$ --- положительные корни уравнения $\varphi_n(x)=1$, при $n\geq c(s)$
их число $p$ постоянно и не больше $\Big[\frac{r}{2}\Big]$. }

Случай малых $s$ и $r$ (при всех $n\in\mathbb{N}$) см. в [\ref{Kotova_Trigub2012}].

Модули гладкости изучали и для аналитических функций на множествах комплексной
плоскости. См. П.~М.~Тамразов [\ref{Tamrazov}] (модуль определяется разделенными
разностями, см. также [\ref{Dzyadyk_Shevchuk}]).

Переходим теперь к рядам Фурье степенного типа, т.е. к функциям из пространства Харди
$H_p(D)$, $D=\big\{z: |z|<1\big\}$, $p>0$.

Введем разные модули гладкости.

$f(e^{it})$ --- предельная функция на окружности $\partial D$. При $r\in\mathbb{N}$ и
$h>0$ контурный или граничный модуль определяется так:
\begin{equation*}
    \omega_r(f;h)_p=\sup\limits_{0<\delta\leq
    h}\Big\|\sum\limits_{\nu=1}^r\binom{r}{\nu}(-1)^\nu
    f\big((\cdot)e^{i\nu\delta}\big)\Big\|_{H_p(D)}.
\end{equation*}

Радиальный модуль $\Big(h\in\Big(0,\frac{2}{p}\Big]\Big)$
\begin{equation*}
    \omega_r(f;rad,h)_p=\Bigg(\int_{-\pi}^\pi\Big|\sum\limits_{\nu=0}^r\binom{r}{\nu}(-1)^\nu f\big(e^{it}(1-\nu
    h)\big)\Big|^pdt\Bigg)^{\frac{1}{p}}.
\end{equation*}

И наконец, линеаризованный граничный модуль $(q\in\mathbb{N})$
\begin{equation*}
    \tilde{\omega}_(f;h)_p=\Bigg\|\int\limits_{[0,1]^q}\sum\limits_{\nu=0}^r\binom{r}{\nu}(-1)^\nu f\Big((\cdot)e^{i\nu
    h\sum\limits_{j=1}^qu_j}\Big)du\Bigg\|_{H_p}
\end{equation*}
($q$-кратный интеграл по кубу).

\textbf{2.8.} \emph{Предположим, что $r\in\mathbb{N}$, $p\in(0,+\infty]$, $f\in H_p(D)$
и $h\in\Big(0,\frac{1}{r+1}\Big]$.}

\emph{а) Если $q=1$ при $p\geq1$ и $q=\Big[\frac{1}{p}+\frac{1}{2}\Big]$ при
$p\in(0,1)$, то
\begin{equation*}
    \omega_r(f;h)_p\asymp\tilde{\omega}_r(f;h)_p.
\end{equation*}
}

\emph{б) Если $S_0(z)\equiv0$, а
$S_{r-1}(z)=\sum\limits_{k=0}^{r-1}\frac{1}{k!}f^{(k)}(0)z^{k}$ при $r\geq2$, то
\begin{equation*}
\omega_r(f;rad,h)_p\asymp\omega_r(f-S_{r-1};h)_p
\end{equation*}
(двойные неравенства с положительными константами, зависящими только от $r$ и $p$).}

См. \textbf{8.4.1} в [\ref{Trigub_Belinsky}]. Там же показано применение новых модулей
гладкости к неравенству Харди-Литтльвуда о росте модулей производных функций при подходе
к границе $\partial D$.

Учитывая, что при $r=1$ $\ \omega_1(f,rad;h)_p=\|f-f_{1-h}\|_p$ (приближение средними
Абеля-Пуассона), получаем при $p\in(0,1]$ и тем более при $p>1$ такое следствие:
\begin{equation*}
\big\|f-f_{1-h}\big\|_{H_p}\asymp\sup\limits_{0<\delta\leq
h}\big\|f-f_{1-\delta}\big\|_{H_p}\asymp\omega(f;h)_p.
\end{equation*}

Вернемся к периодическим функциям с полным рядом Фурье.

\textbf{2.9.} \emph{а) Всегда
$\lim\limits_{h\rightarrow0}\frac{\omega_r(f;h)_p}{h^r}=\sup\limits_{h\in\big[0,\frac{1}{r}\big]}\frac{\omega_r(f;h)_p}{h^r}$.
А для того чтобы $\omega_r(f;h)_p\leq h^r$ или, что то же самое,
$\tilde{\omega}_r(f;h)\leq\frac{h^r}{r+1}$ необходимо и достаточно, чтобы существовала
абсолютно непрерывная производная $f^{(r-1)}$ и $\big\|f^{(r)}\big\|_p\leq1$ при
$p\in(1,+\infty]$, а при $p=1$ полная вариация $f^{(r-1)}$ была не больше единицы }(см.
[\ref{Trigub_Belinsky}]).

\emph{б) Для любой функции из $L_p$, $p\in(0,1)$
\begin{equation*}
    \lim\limits_{h\rightarrow0}\frac{\omega_r(f;h)_p}{h^{r-1+\frac{1}{p}}}=\sup\limits_{0<h\leq\frac{\pi}{2}}\frac{\omega_r(f;h)_p}{h^{r-1+\frac{1}{p}}}.
\end{equation*}
}

\emph{А $\omega_r(f;h)_p\leq h^{r-1+\frac{1}{p}}$ тогда и только тогда, когда
$f^{(r-1)}$ --- ступенчатая функция или кусочно постоянная (возможно, после исправления
на множестве нулевой меры) и при любом наборе точек $x_0<x_1<...<x_n$
\begin{equation*}
    \sum\limits_{k=0}^{n-1}\Big|f^{(r-1)}(x_k)-f^{(r-1)}(x_{k+1})\Big|^p\leq\Big(\int\limits_0^r\frac{1}{(r-1)!}\big|\Delta_1^r(-t)_+^{r-1}\big|^pdt\Big)^{-1}.
\end{equation*}
}

См. [\ref{Kolomoitsev}].

\textbf{Замечание 1} (о других модулях)

Иногда рассматривают на прямой следующее условие
\begin{equation*}
    {\rm Lip}\ \alpha:\quad \big|f(x)-f(y)\big|\leq\frac{|x-y|^\alpha}{\big(1+|x|\big)^\alpha\big(1+|y|\big)^\alpha}
\end{equation*}
(см. следствие из \textbf{5.1} в [\ref{Liflyand_Samko_Trigub}]).

А для функций на отрезке $[-1,1]$, напр., крайние точки играют особую роль:
\begin{equation*}
\big|f(x)-f(y)\big|\leq\Bigg(\frac{|x-y|}{\sqrt{1-x^2}+\sqrt{1-y^2}+|x-y|}\Bigg)^\alpha
\end{equation*}
(см. [\ref{Trigub_Belinsky}], с.~172).

Если $\varphi(x)=\sqrt{1-x^2}$, то Д-Т модуль функции $f\in L_p[-1,1]$,
$p\in[1,+\infty]$, равен
\begin{equation*}
    \omega_r^\varphi(f;h)_p=\sup\limits_{0<\delta\leq
    h}\big\|\dot{\Delta}_{\delta\varphi}^rf(\cdot)\big\|_p,
\end{equation*}
где $L_p$-норма по отрезку $[-1,1]$ в предположении, что
$\dot{\Delta}_{\delta\varphi}^r=0$, если $r\varphi\delta$ или $-r\varphi\delta$ не
принадлежат $[-1,1]$,
\begin{equation*}
\inf\limits_g\Big(\|f-g\|_p+h^r\big\|\varphi^rg^{(r)}\big\|_p\asymp\omega_r^\varphi(f;h)_p
\end{equation*}
(см. [\ref{Ditzian_Totik}], а также [\ref{De_Vore}]).

\begin{center}
    \textbf{\S3 Модули гладкости и $K$-функционалы в $\mathbb{R}^d$ ($d\geq2$)}
\end{center}

Для функции $f\in L_p$, $p\in[1,+\infty]$ модуль гладкости естественно ввести следующим
образом ($E\subset\mathbb{R}^d$, $E\neq\varnothing$)
\begin{equation*}
    \omega_r(f;E;h)_p=\sup\limits_{u\in
    E}\Big\|\sum\limits_{\nu=0}^r\binom{r}{\nu}(-1)^\nu f(\cdot+\nu hu)\Big\|_p.
\end{equation*}

Для монотонности по шагу $h>0$ множество $E$ считаем звездным относительно нуля. Будем
предполагать еще, что $E$ --- компакт в шаре радиуса 1.

Если $E$ --- единичный шар, то будем писать $\omega_r^\circ$ (наибольший модуль).
Наименьший модуль получается, если в качестве $E$ взять отрезок, выходящий из нуля
(модуль по направлению). Если $E$ --- объединение $d$ единичных отрезков в направлении
прямоугольных осей координат, то будем писать $\omega_r^+(f;h)$. Именно с таким модулем
формулируют обычно прямую теорему типа Джексона (см. [\ref{De_Vore}],
[\ref{Trigub_Belinsky}]).

Определение $\Phi_\varepsilon(f)$ см. в конце $\S 1$. 3.1--3.4 есть в
[\ref{Trigub1980}], [\ref{Trigub1988}] (см. также [\ref{Trigub_Belinsky}], гл.~VIII).

\textbf{3.1.} \emph{Неравенство
\begin{equation*}
    c_1(r,d)\omega_r^+(f;\varepsilon)\leq\big\|f-\Phi_\varepsilon(f)\big\|_\infty\leq
    c_2(r,d)\omega_r^\circ(f;\varepsilon)
\end{equation*}
не может быть верным для всех функций из $C(\mathbb{T}^d)$ при $r\in\mathbb{N}$ и
$d\geq2$, если $\varphi$ непрерывна на своем носителе,содержащемся в кубе $[-1,1]^d$.  }

Введем теперь линеаризованный модуль гладкости четного порядка ($\mu$ --- конечная
комплекснозначная борелевская мера)
\begin{equation*}
    \tilde{\omega}_{2r}(f;\mu;h)=\Big\|\int\limits_{\mathbb{R}^d}\sum\limits_{\nu=0}^{2r}\binom{2r}{\nu}(-1)^\nu
    f\big(\cdot+(\nu-r)hu\big)d\mu(u)\Big\|
\end{equation*}
(интеграл от симметричной разности $\dot{\Delta}^{2r}$ по мере $\mu$).

Если $d\mu=\chi_E(u) du$ ($\chi_E$ --- индикатор $E$), то в случае, когда $E$ ---
единичный шар, пишем $\tilde{\omega}_{2r}^\circ(f;h)$, а в случае
$d\mu=\sum\limits_{j=1}^d\chi_{E_j}du_j$,  $E_j=[-1,1]$ на оси $Ox_j$,\  $1\leq
    j\leq d,$ пишем $\tilde{\omega}_{2r}^+(f;h)$ ($L_\infty=C$).

\textbf{3.2.} \emph{Пусть множество $E$ обладает следующей симметрией: любая точка $E$
остается в $E$ после перестановки местами двух координат или перемены знака любой
координаты. Тогда при $p\in[1,+\infty]$
\begin{equation*}
    \tilde{\omega}_2(f,E,h)_p\asymp\tilde{\omega}^\circ_2(f,E,h)_p.
\end{equation*}
}

\textbf{3.3.} \emph{При $r\geq2$ и $d\geq2$ модули $\tilde{\omega}^\circ_{2r}$ и
$\tilde{\omega}^\Box_{2r}$ ($E$ --- единичный квадрат) не сравнимы в $C(\mathbb{T}^d)$
при $h\rightarrow0$.}

\textbf{3.4.} \emph{Если $\Delta$ --- оператор Лапласа, $r\in\mathbb{N}$ и
$\delta>\frac{d-1}{2}$, то при $p\in[1,+\infty]$ и $\varepsilon>0$
\begin{equation*}
    K\big(\varepsilon^{2r},f,L_p,\Delta^r\big)_p\asymp\Big\|f-\sum\limits_{k\in\mathbb{Z}^d}\big(1-\varepsilon^{2r}|k|^{2r}\big)_+^\delta\hat{f}_ke_k\Big\|_p\asymp\tilde{\omega}_{2r}^\circ(f,\varepsilon)_p.
\end{equation*}
}

\textbf{3.5.} \emph{Если $D_{2r}=\sum\limits_{j=1}^d\frac{\partial^{2r}}{\partial
x_j^{2r}}$, $r\in\mathbb{N}$ и $\delta>\frac{d-1}{2}$, то при $p\in [1,+\infty]$ и
$\varepsilon>0$
\begin{equation*}
K\big(\varepsilon^{2r},f,L_p,D_{2r}\big)_p\asymp\Big\|f-\sum\limits_{k\in\mathbb{Z}^d}\Big(1-\varepsilon^{2r}\sum\limits_{j=1}^dk_j^{2r}\Big)_+^\delta\hat{f}_ke_k\Big\|_p\asymp\tilde{\omega}_{2r}^+(f,\varepsilon)_p.
\end{equation*}
}

\textbf{3.6.} \emph{Если $S_\nu^\square(f)$, $\nu\in\mathbb{Z}_+$ --- $\nu$-тая
квадратная частная сумма ряда Фурье функции $f\in C(\mathbb{T}^2)$, а
\begin{equation*}
    d_1(f)\sim\sum\limits_{k\in\mathbb{Z}^2}\max\big\{|k_1|,|k_2|\big\}\hat{f}_ke_k,
\end{equation*}
то для приближения суммами Марцинкевича
\begin{equation*}
    \begin{split}
        & \Big\|f-\frac{1}{n+1}\sum\limits_{\nu=0}^nS_\nu^\square(f)\Big\|_\infty\asymp\Big\|\int\limits_1^\infty\Big(\dot{\Delta}_{t(e_1^\circ+e_2^\circ)/n}^2+\dot{\Delta}_{t(e_1^\circ-e_2^\circ)/n}^2\Big)f(\cdot)\frac{dt}{t^2}\Big\|_\infty\asymp \\
        &\asymp K\Big(\frac{1}{n};f,C,W(d_1)_\infty\Big),
    \end{split}
\end{equation*}
где $(e_1^\circ,e_2^\circ)$ --- стандартный базис в $\mathbb{R}^2$, а под интегралом ---
сумма вторых симметричных разностей в направлении биссектрис первой и четвертой четверти
$\mathbb{R}^2$.
 }

 Первое соотношение доказано О.~И.~Кузнецовой (см. [\ref{Kuznetsova_Trigub}]).

\textbf{3.7.} \emph{Если
$\Delta_{r,\delta}^+f(x)=\sum\limits_{j=1}^d\sum\limits_{\nu=0}^{2r}\binom{2r}{\nu}(-1)^\nu
f\big(x+(\nu-r)\delta e_j^\circ\big)$, где $\{e_j^\circ\}_1^d$ --- стандартный базис, то
при $r\in \mathbb{N}$ и $p\in[1,+\infty]$
\begin{equation*}
\tilde{\omega}_{2r}^+(f,h)_p\asymp\sup_{0<\delta\leq h}\big\|\Delta_{r,\delta}^+
f(\cdot)\big\|_p,
\end{equation*}
\begin{equation*}
\tilde{\omega}_{2r}^\circ(f,h)_p\asymp\sup_{0<\delta\leq h}\big\|(\Delta_{1,\delta}^+)^r
f(\cdot)\big\|_p
\end{equation*}
}

\textbf{3.8.} \emph{При любом $\alpha>0$, $\beta>\frac{d-1}{2}$, натуральном
$r>\frac{\alpha}{2}$, $p\in[1,+\infty]$ и
\begin{equation*}
    d_\alpha(f)\sim\sum\limits_{k\in\mathbb{Z}^d}\sum\limits_{j=1}^d|k_j|^\alpha\hat{f}_ke_k
\end{equation*}
\begin{equation*}
    \begin{split}
        & K\big(\varepsilon^\alpha,f,L_p,W(d_\alpha)_p\big)\asymp\Big\|f-\sum\limits_{k\in\mathbb{Z}^d}\Big(1-\sum\limits_{j=1}^d\varepsilon^\alpha|k_j|^\alpha\Big)^\beta_+\hat{f}_ke_k\Big\|_p\asymp \\
        & \asymp\Big\|\sum\limits_{j=1}^d\int\limits_1^\infty\frac{1}{u^{1+\alpha}}\dot{\Delta}_{\varepsilon ue_j}^{2r}f(\cdot)du\Big\|_p
    \end{split}
\end{equation*}
(двусторонние неравенства с положительными константами, зависящими только от $d$,
$\alpha$, $\beta$ и $r$) }

\textbf{3.9.} \emph{При любом $\alpha>0$, $\beta>\frac{d-1}{2}$, натуральном
$r>\frac{1}{2}(\alpha+d-1)$ и
\begin{equation*}
    d_\alpha(f)\sim\sum\limits_{k\in\mathbb{Z}^d}(\varepsilon|k|)^\alpha\hat{f}_ke_k
\end{equation*}
для всех $f\in L_p(\mathbb{T}^d)$ и $\varepsilon>0$
\begin{equation*}
    \begin{split}
        & K\big(\varepsilon^\alpha,f,L_p,W(d_\alpha)_p\big)\asymp\Big\|f-\sum\limits_{k\in\mathbb{Z}^d}\Big(1-\varepsilon^\alpha|k_j|^\alpha\Big)_+\hat{f}_ke_k\Big\|_p\asymp \\
        & \asymp\Big\|\int\limits_{|u|\geq1}\frac{1}{|u|^{\alpha+r}}\dot{\Delta}_{\varepsilon
        u}^{2r}f(\cdot)du\Big\|_p.
    \end{split}
\end{equation*}}

Теоремы 3.8 и 3.9 анонсированы автором еще в статье [\ref{Kuznetsova_Trigub}], но
доказательство приводится впервые. Так же доказываются и другие теоремы, приведенные в
настоящей статье.

\textbf{Доказательство теоремы 3.9.}

$\rhd$ Сначала, применяя принцип сравнения (см. в конце $\S 1$), докажем второе
соотношение (о точном порядке приближения), а затем, используя лемму (см. там же),
получим первое соотношение.

Понадобятся некоторые свойства винеровских банаховых алгебр $A(\mathbb{R}^d)$ и
$B(\mathbb{R}^d)$.

Очевидно, что функции из $B$ ограниченные и равномерно непрерывные, а если $f\in A$
(идеал в $B$), то $f\in C_0(\mathbb{R}^d)$
($f(\infty)=\lim\limits_{|x|\rightarrow\infty}f(x)=0$ в силу леммы Римана-Лебега).

Алгебра $B(\mathbb{R}^d)$ с ростом $d$ расширяется.

\textbf{А.} Начнем с необходимых условий.

Если $f\in A(\mathbb{R}^d)$, то ее радиальная часть $f_0(t)$, $t\in (0,+\infty)$
(интегральное среднее по сфере $|x|=t$) удовлетворяет следующим условиям: $f_0\in
C^{d_1}(0,+\infty)$, $d_1=\Big[\frac{d-1}{2}\Big]$,
\begin{equation*}
\lim\limits_{t\rightarrow+\infty}t^kf_0^{(k)}(t)=0\qquad (0\leq k\leq d_1),
\end{equation*}
а при $d\geq3$ и
\begin{equation*}
\lim\limits_{t\rightarrow+0}t^kf_0^{(k)}(t)=0\qquad (0\leq k\leq d_1).
\end{equation*}

Кроме того, при любом $t>0$ сходится интеграл
\begin{equation*}
    \int\limits_{\rightarrow0}^t\frac{f_0^{(d_1)}(t+u)-f_0^{(d_1)}(t-u)}{u^{\frac{d+1}{2}-d_1}}du
\end{equation*}
(не обязательно абсолютно). См. [\ref{Trigub1980}], \textbf{6.5.7} в
[\ref{Trigub_Belinsky}].

\textbf{Б.} Функция $f$ принадлежит $A(B)$ в точке, если она допускает продолжение с
некоторой ее окрестности до функции из $A(\mathbb{R}^d)$ ($B(\mathbb{R}^d)$). Как
известно, если $f$ принадлежит $A(B)$ в любой точке, включая $\infty$ (окрестность
$\infty$ --- это $|x|>M$), то $f\in A(\mathbb{R}^d)$ ($B(\mathbb{R}^d)$) (локальные
свойства). Разница между функциями из $A$ и $B$ только в окрестности $\infty$.

Если $f\in B(\mathbb{R}^d)$, то $\frac{1}{|x|^\alpha}f(x)\in A(\mathbb{R}^d)$ при
$\alpha>0$ (см. ниже в самом конце в \textbf{Д}).

Используется ниже и $\frac{1}{f}$-теорема Винера (см., напр.,
[\ref{Liflyand_Samko_Trigub}]): если $f(x)-f(\infty)\in A (\mathbb{R}^d)$,
$f(\infty)\neq0$ и $f(x)\neq0$ для $x\in \mathbb{R}^d$, то и
$\frac{1}{f(x)}-\frac{1}{f(\infty)}\in A(\mathbb{R}^d)$.

\textbf{В.} Достаточные условия принадлежности $A(\mathbb{R}^d)$.

Если $y=(y_1,...,y_d)$ и $y_j\neq0$ ($1\leq j\leq d$), то при любом $h\in \mathbb{R}^d$
\begin{equation*}
    \big\|f\big((y,\cdot)+h\big)\big\|_B=\|f\|_B.
\end{equation*}

Если ряд Фурье непрерывной периодической функции $f$ сходится абсолютно, то $f\in
B(\mathbb{R}^d)$ и $\|f\|_B=\sum\limits_{k\in \mathbb{Z}^d}\big|\hat{f}_k\big|$.

А для этого достаточно, чтобы гладкость функции в $L_2(\mathbb{T}^d)$ была больше
$\frac{d}{2}$ (см., напр., [\ref{Stein_Weiss}], п.~3, гл.~VII).

Положим при $0<a<b$ $\ \chi_{a,b}(x)=1$ при $|x|\leq a$, $\chi_{a,b}(x)=0$ при $|x|\geq
b$ и $\chi_{a,b}\in C^\infty(\mathbb{R}^d)$.

В силу формулы обращения преобразования Фурье при любом $d\geq1$ $\ \chi_{a,b}\in
A(\mathbb{R}^d)$.

Если в некоторой окрестности точки $x_0\in \mathbb{R}^d$ функция $f$ имеет гладкость в
$L_2$ (еще лучше --- в $C$) больше $\frac{d}{2}$, то она принадлежит $A(\mathbb{R}^d)$ в
этой окрестности.

Действительно, после умножения на $\chi_{a,b}(x-x_0)$ при достаточно малых $a$ и $b$ она
допускает продолжение до периодической функции из $B(\mathbb{R}^d)$ и после умножения
продолженной функции на $\chi_{a,b}(x-x_0)$ получаем функцию из $A(\mathbb{R}^d)$.

\textbf{Г.} Дополнение к \textbf{В}.

При $d=1$, как заметил А.~Зигмунд, можно предполагать для абсолютной сходимости ряда
Фурье гладкость в $C$ только больше нуля, если добавить гладкость в $L_1(\mathbb{T})$
(ограниченность вариации).

Если при $d\geq2$ добавить выпуклость, то гладкость в $C$ больше $\frac{d}{2}$ можно
заменить на гладкость в $C$ больше $\frac{d-1}{2}$; точнее, если $f$, как функция от
$x_j$ ($1\leq j\leq d$), имеет в окрестности точки частую производную
$\frac{\partial^{d_1}}{\partial x^{d_1}}f$ $\Big(d_1=\Big[\frac{d-1}{2}\Big]\Big)$,
удовлетворяющую условию Липшица степени больше $\frac{d+1}{2}-d_1-1$ равномерно по
остальным $(d-1)$ переменным и выпуклую по $x_j$ (вверх или вниз) или меняющую
выпуклость при переходе через точку, то $f$ в этой точке принадлежит $A$ (см. более
точную теорему в [\ref{Trigub1980}] \textbf{4.III} или в [\ref{Trigub_Belinsky}],
\textbf{6.4.5}).

\textbf{Д.} Радиальные функции $f_0\big(|x|\big)$.

Для того чтобы $F_0\big(|x|\big)\in A(\mathbb{R}^d)$ при $d\geq2$ необходимо и
достаточно, чтобы существовала функция $f_0\big(|x|\big)\in A(\mathbb{R}^1)$ такая, что
при $t\geq0$
\begin{equation*}
    F_0(t)=\int\limits_0^1f_0(ut)\big(1-u^2\big)^{\frac{d-3}{2}}du.
\end{equation*}

Функция $f_0$ по $F_0$ определяется однозначно: при нечетном $d$ простым
дифференцированием, а при четном $d$ --- дифференцированием полуцелого порядка (см.
[\ref{Trigub_Belinsky}], \textbf{6.3.6} и [\ref{Trigub2010}], соответственно). См. также
\textbf{6.5.8} в [\ref{Trigub_Belinsky}] (используется и асимптотика преобразования
Фурье радиальных функций с условиями типа выпуклости).

Преобразование Фурье радиальной функции есть функция радиальная, а множество радиальных
функций из $A(\mathbb{R}^d)$ с ростом $d$ уменьшается.

При нечетном $d$ из \textbf{Б$_1$} ($\S2$) можно вывести следующее достаточное условие.

Если $f_0\in C[0,+\infty)$, $f_0^{(d_1)}\in AC_{loc}(0,+\infty)$ и
$\int\limits_0^\infty\underset{u\geq t}{{\rm
ess}\sup}\Big|u^{d_1}f_0^{(d_1+1)}(u)\Big|dt<\infty$, то $f_0\big(|x|\big)\in
A(\mathbb{R}^{2d_1+1})$ тогда и только тогда, когда
\begin{equation*}
    \lim\limits_{t\rightarrow+\infty}t^kf_0^{(k)}(t)=0\qquad (0\leq k\leq
    d_1=\frac{d-1}{2}).
\end{equation*}

Из приведенной выше формулы видно, что если $f_0(t)=t^\gamma$ и $\gamma>-1$, то и
$F_0(t)=ct^\gamma$, $c\neq0$. Применяем это в окрестности нуля и $\infty$.

\textbf{Пример.} $\big(1-|x|^\alpha\big)_+^\beta\in A(\mathbb{R}^d)$ при $\alpha>0$ и
$\beta>\frac{d-1}{2}$ и только в этом случае.

Необходимость условия на $\beta$ следует из \textbf{А} при $t=1$.

При любом $m\in \mathbb{N}$ в достаточно малой окрестности нуля
$\big(1-|t|^\alpha\big)_+^\beta=\sum\limits_{k=0}^ma_{k,\alpha,\beta}|t|^{k\alpha}+h(|t|)$
$\big(h\in C^n$ при $n<\alpha(m+1)\big)$. И выбираем $m$ так, чтобы
$\alpha(m+1)>\frac{d}{2}+1$ (см. \textbf{В}). В точках, в которых $|x|=1$, применяем
\textbf{Г}.

Переходим теперь непосредственно к доказательству второго соотношения в теореме~3.9.

В силу принципа сравнения (см. в конце $\S\ 1$) нужно проверить, что
$g=\frac{\psi}{\varphi}$ и $\frac{1}{g}\in B(\mathbb{R}^d)$, где
\begin{equation*}
    \psi(x)=\int\limits_{|u|\geq1}\frac{\sin^{2r}(x,u)}{|u|^{\alpha+d}}du,\qquad
    \varphi(x)=1-\big(1-|x|^\alpha\big)_+^\beta.
\end{equation*}

Очевидно, что при $x\neq0$ $\ \psi(x)$ и $\varphi(x)>0$. А то, что
$g(0)=\lim\limits_{x\rightarrow0}g(x)\neq0$, будет доказано позже.

Применяя к $\sin$ формулу Эйлера и бином Ньютона, получим
\begin{equation*}
    \psi(x)=\frac{1}{\alpha}\binom{2r}{r}+\sum\limits_{\nu=0,\ \nu\neq
    r}^{2r}\binom{2r}{\nu}(-1)^{\nu+r}\int\limits_{|u|\geq1}\frac{1}{|u|^{\alpha+d}}e^{2i(r-\nu)(x,u)}du.
\end{equation*}

Отсюда $\psi(\infty)=\frac{1}{\alpha}\binom{2r}{r}$, а $\psi-\psi(\infty)\in
A(\mathbb{R}^d)$ (по определению алгебры $A$).

Но (см., напр., [\ref{Stein_Weiss}])
\begin{equation*}
\int\limits_{|u|\geq1}\frac{1}{|u|^{\alpha+d}}e^{i(x,u)}du=\frac{c_1(d)}{|x|^{\frac{d}{2}-1}}\int\limits_{1}^\infty\frac{1}{|u|^{\alpha+d}}u^{\frac{d}{2}}J_{\frac{d}{2}-1}\big(u|x|\big)du,
\end{equation*}
где $J_\mu$ --- функция Бесселя первого рода. Как известно,
\begin{equation*}
    J_\mu\in C^\infty(\mathbb{R}\setminus\{0\}),\qquad |J_\mu(t)|\leq
    c_2(\mu)\min\Big\{|t|^\mu,\frac{1}{\sqrt{|t|}}\Big\}.
\end{equation*}

Поэтому $\psi\in C^\infty(\mathbb{R}^d\setminus\{0\})$, так как указанный интеграл равен
\begin{equation*}
    c_1(d)|x|^{\alpha-1}\int\limits_{|x|}^\infty\frac{1}{u^{\alpha+\frac{d}{2}}}J_{\frac{d}{2}-1}(u)du.
\end{equation*}

С другой стороны, при $\alpha>0$ и $\beta>\frac{d-1}{2}$ (см. пример в \textbf{Д})
\begin{equation*}
    \big(1-|x|^\alpha\big)_+^\beta\in A(\mathbb{R}^d).
\end{equation*}

Но тогда по $\frac{1}{f}$-теорема Винера (см. \textbf{Б}) при $|x|>\delta>0$
\begin{equation*}
    \frac{1}{\varphi(x)}-1=\frac{1}{1-\big(1-|x|^\alpha\big)_+^\beta}-1\in
    A(\mathbb{R}^d).
\end{equation*}

Так что при $|x|\geq\delta>0$ $\ g-g(\infty)=\frac{\psi}{\varphi}-\psi(\infty)\in
A(\mathbb{R}^d)$.

Остается окрестность нуля.
\begin{equation*}
    \psi(x)=\int\limits_{\mathbb{R}^d}\frac{\sin^{2r}(x,u)}{|u|^{\alpha+d}}du-\int\limits_{|u|\leq1}\frac{\sin^{2r}(x,u)}{|u|^{\alpha+d}}du.
\end{equation*}

Первый интеграл после замены переменных (вращение $\mathbb{R}^d$) равен
\begin{equation*}
    \int\limits_{\mathbb{R}^d}\frac{\sin^{2r}(x,u)}{|u|^{\alpha+d}}du=\int\limits_{\mathbb{R}^d}\frac{\sin^{2r}(|x|\cdot|u|)}{|u|^{\alpha+d}}du=\int\limits_{0}^\infty\frac{\sin^{2r}u|x|}{u^{\alpha+1}}du=|x|^\alpha\int\limits_{0}^\infty\frac{\sin^{2r}u}{u^{\alpha+1}}du.
\end{equation*}

Во втором интеграле разлагаем $\sin^{2r}t$ в степенной ряд, учитывая, что
\begin{equation*}
    \int\limits_{|u|\leq1}(x,u)^{2k}du=|x|^{2k}\frac{(2k-1)!!}{2^k}\frac{\pi^{\frac{k}{2}}}{\Gamma\Big(k-1+\frac{d}{2}\Big)}\qquad
    (k\in\mathbb{N})
\end{equation*}

Получаем после вращения $\mathbb{R}^d$ (по условию $2r>\alpha+d-1$)
\begin{equation*}
    \begin{split}
        & \int\limits_{|u|\leq1}\frac{\sin^{2r}(x,u)}{|u|^{\alpha+d}}du=\sum\limits_{k=r}^\infty c_{k,r}\int\limits_{|u|\leq1}\frac{(x,u)^{2k}}{|u|^{\alpha+d}}du=\sum\limits_{k=r}^\infty
        c_{k,r}|x|^{2k}\int\limits_{|u|\leq1}\frac{|u|^{2k}}{|u|^{\alpha+d}}du= \\
        &=\sum\limits_{k=r}^\infty
        c_{k,r}|x|^{2k}\int\limits_0^1u^{2k-\alpha-1}du=\sum\limits_{k=r}^\infty
        c_{k,r}\frac{1}{2k-\alpha}|x|^{2k}.
    \end{split}
\end{equation*}

Так что
\begin{equation*}
    \frac{\psi(x)}{|x|^\alpha}=\int_0^\infty\frac{\sin^{2r}u}{u^{1+\alpha}}du+|x|^{2r-\alpha}\psi_1(|x|),
\end{equation*}
где $2r>\alpha$, а $\psi_1$ --- целая функция.

Следовательно, в окрестности нуля (см. еще степенной пример $F_0$ в \textbf{Д})
\begin{equation*}
    \frac{\psi(x)}{|x|^\alpha}\in A(\mathbb{R}^d).
\end{equation*}

Но и
\begin{equation*}
    \frac{\varphi(x)}{|x|^\alpha}=\frac{1-\big(1-|x|^\alpha\big)_+^\beta}{|x|^\alpha}=\varphi_1(|x|^\alpha),
\end{equation*}
где $\varphi_1(0)=\beta$ и $\varphi_1$ --- аналитична в окрестности нуля.

Следовательно, и $\frac{1}{\varphi_1(|x|^\alpha)}\in A(\mathbb{R}^d)$ в малой
окрестности нуля.

Поэтому и в окрестности нуля
\begin{equation*}
    g(x)=\frac{\psi(x)}{\varphi(x)}=\frac{\psi(x)}{|x|^\alpha}\cdot\frac{|x|^\alpha}{\varphi(x)}\in
    A(\mathbb{R}^d).
\end{equation*}

Следовательно, и в целом
\begin{equation*}
    g(x)-g(\infty)=\frac{\psi(x)}{\varphi(x)}-\psi(\infty)\in A(\mathbb{R}^d).
\end{equation*}

Но тогда по теореме Винера и $\frac{1}{g}-\frac{1}{\psi(\infty)}\in A(\mathbb{R}^d)$.

Второе соотношение в 3.9 доказано.

Для вывода формулы для $K$-функционала понадобится лемма, сформулированная в конце
$\S1$.

\textbf{Доказательство леммы.}

$\rhd$ Оценка сверху.

В силу определения $K$-функционала
\begin{equation*}
    K\big(\varepsilon,f,L_p,W(d_\alpha)_p\big)\leq
    \|f-\Phi_\varepsilon(f)\|_p+\varepsilon\|d_\alpha(\Phi_\varepsilon(f))\|_p.
\end{equation*}

Применяем условие $\gamma)$.

Оценка снизу.
\begin{equation*}
\|f-\Phi_\varepsilon(f)\|_p=\big\|(f-g)-\Phi_\varepsilon(f-g)+\big(g-\Phi_\varepsilon(g)\big)\big\|_p.
\end{equation*}

Применяем условия $\beta)$ и $\alpha)$:
\begin{equation*}
\|f-\Phi_\varepsilon(f)\|_p\leq(1+a)\|f-g\|+a\varepsilon\|d_\alpha
g\|\leq(1+a)\big(\|f-g\|+\varepsilon\|d_\alpha g\|\big).
\end{equation*}

Осталось взять нижнюю грань по $g\in W(d_\alpha)$.
$\qquad\qquad\qquad\qquad\qquad\qquad\qquad\qquad\blacktriangleleft$

Заметим, что вид оператора $d_\alpha$ по $\Phi_\varepsilon$ определяется по классу
насыщения для $f-\Phi_\varepsilon(f)$.

Теперь докажем первое соотношение в 3.9.

Для доказательства $\alpha)$ в силу принципа сравнения нужно проверить, что
$\frac{\varphi(x)}{|x|^\alpha}\in A(\mathbb{R}^d)$. В окрестности нуля это уже есть. А
то, что $\frac{1}{|x|^\alpha}\in A(\mathbb{R}^d)$ при $|x|\geq\delta>0$ и любом $d$,
давно известно (см. [\ref{Liflyand_Samko_Trigub}]). Можно также при нечетном $d$ (см.
\textbf{Е}) перейти к случаю $d=1$ (будет функция того же вида) и применить
\textbf{B$_1$} ($\S2$).

Условия $\beta)$ выполняется всегда, если есть сходимость в $C(\mathbb{T}^d)$ (по
теореме Банаха--Штейнгауза).

Условие $\gamma)$ выполняется, так как уже доказано, что
\begin{equation*}
    \frac{|x|^\alpha\big(1-|x|^\alpha\big)_+^\beta}{\varphi(x)}\in A(\mathbb{R}^d).
\end{equation*}

Теорема 3.9 полностью доказана.

\textbf{Замечание 2} (о случае $L_p$, $p\in(1,+\infty)$).

В случае пространства $C(\mathbb{T}^d)$ и компактности оператора приближения справа в
принципе сравнения условие $\|g\|_B<\infty$ для переходной функции не только достаточно,
но и необходимо (см. [\ref{Trigub1980}] и [\ref{Trigub_Belinsky}]). Так что из такого
неравенства в $C(\mathbb{T}^d)$ следует сразу такое же неравенство и для всех функций из
$L_p(\mathbb{T}^d)$, $p\geq1$. Но в случае пространства $L_p$, $p\in(1,+\infty)$, есть
более слабые достаточные условия (напр., теорема Марцинкевича о мультипликаторах, см.
[\ref{Stein}]).

Таким образом получаются, напр., следующие неравенства
\begin{equation*}
    K\big(\varepsilon^{2r};f;L_p,\Delta^r\big)\asymp
    \tilde{\omega}_{2r}^\circ(f,\varepsilon)_p\asymp\omega_{2r}^\circ(f,\varepsilon)_p.
\end{equation*}

(См. \textbf{8.2.9} и \textbf{8.3.2 в)} в [\ref{Trigub_Belinsky}]).

\textbf{Замечание 3} (о непериодических функциях).

Модули гладкости функции $f:\mathbb{R}^d\rightarrow\mathbb{C}$ из $L_p(\mathbb{R}^d)$
определяются обычным образом, как и $K$-функционалы. Мультипликатор интегралов Фурье в
$L_p$ как ограниченный линейный оператор, определяемый на функциях из $L_2\cap L_p$
одной функцией $\varphi:\mathbb{R}^d\rightarrow\mathbb{C}$ (см. [\ref{Stein}], гл.~IV),
которую, как доказано в [\ref{Lebedev_Olevskii}], можно, не уменьшая общности, считать
ограниченной и непрерывной почти всюду, следующим образом:
\begin{equation*}
    f(x)\mapsto\int\limits_{\mathbb{R}^d}\varphi(y)\hat{f}(y)e^{i(x,y)}dy.
\end{equation*}

А метод $\Phi_\varepsilon(f)$ суммирования интегралов Фурье равен
\begin{equation*}
    \int\limits_{\mathbb{R}^d}\varphi(\varepsilon y)\hat{f}(y)e^{i(x,y)}dy\qquad
    \Big(\varphi(0)=\frac{1}{(2\pi)^d}\Big).
\end{equation*}

В принципе сравнения стоит норма мультипликатора для переходной функции. А для норм
мультипликаторов рядов и интегралов Фурье, определяемых одной и той же функцией
$\varphi$, есть простые связи в одну и другую стороны (см. [\ref{Stein_Weiss}], гл.~VII,
\textbf{3.8} и \textbf{3.18}).

Подобные приведенным выше соотношения имеют место и для интегралов Фурье.

Приведем лишь один пример [\ref{Kotova_Trigub2015}]:

Если
\begin{equation*}
    \Phi_\varepsilon(f;x)=\frac{1}{(2\pi)^d}\int\limits_{\mathbb{R}^d}f(x+\varepsilon
    y)\hat{\varphi}(y)dy,
\end{equation*}
где
\begin{equation*}
    \varphi(x)=\big(1-|x|^{2r}\big)_+^\delta\qquad \Big(r\in \mathbb{N},
    \delta>\frac{d-1}{2}\Big),
\end{equation*}
то при $p\geq1$ в $L_p(\mathbb{R}^d)$ $(L_\infty=C)$
\begin{equation*}
    \|f-\Phi_\varepsilon(f)\|_p\asymp\Big\|\int\limits_{|u|\leq1}\dot{\Delta}_{\varepsilon u}^{2r}f(\cdot)du\Big\|_p\asymp\inf\limits_g\Big\{\|f-g\|_p+\varepsilon^{2r}\|\Delta^r g\|_p\Big\}.
\end{equation*}

$\Phi_\varepsilon(f)$ --- это средние Бохнера-Рисса.

\begin{center}
$\textbf{\S4}$ \textbf{Случай банахова пространства.}
\end{center}

Пусть $E_1$ и $E_2$ --- два нормированных пространства с нормами $|\cdot|_1$ и
$|\cdot|_2$, соответственно, а $f:E_1\rightarrow E_2$ --- ограниченная функция
$\big(\|f\|_\infty=\sup\limits_{x\in E}|f(x)|_2<\infty\big)$.

Предположим еще, что при любых $x$ и $y\in E_1$
\begin{equation*}
    |f(x+y)|_2\leq|f(x)|_2.
\end{equation*}

Модуль гладкости порядка $r\in\mathbb{N}$ и шага $|h|_1$ определятся обычным образом:
\begin{equation*}
    \omega_r(f,|h|_1)=\sup\limits_{x\in E_1,
    |\delta|_1\leq|h|_1}\big|\Delta_\delta^rf(x)\big|_2, \qquad
    \Delta_\delta^rf(x)=\sum\limits_{\nu=0}^r\binom{r}{\nu}(-1)^\nu f(x+\nu\delta).
\end{equation*}

$\lim\limits_{h\rightarrow0}\omega_r(f,|h|_1)=0$ --- это критерий равномерной
непрерывности $f$ (при $r\geq2$ можно применить обобщенное неравенство Маршо (см. ниже
перед 4.4) $r-1$ раз).

Всегда (см. ниже 4.2)
\begin{equation*}
\lim\limits_{h\rightarrow0}\frac{\omega_r(f,|h|_1)}{|h|_1^r}=\sup\limits_{h\in
E_1}\frac{\omega_r(f,|h|_1)}{|h|_1^r}.
\end{equation*}

Если верхняя грань конечна, то будем писать $f\in W^r$. $W^r\subset W^{r-1}$ (см. ниже
обобщенное неравенство Маршо).

Если $f$ дифференцируема на $E_1$, т.е. существует линейный ограниченный оператор
$f'(x): E_1\rightarrow E_2$ такой, что при $x\in E_1$
\begin{equation*}
\lim\limits_{h\rightarrow0}\Big(\frac{f(x+h)-f(x)}{|h|_1}-f'(x)h^\circ\Big)=0\qquad
\Big(h^\circ=\frac{h}{|h|_1}\Big),
\end{equation*}
то
\begin{equation*}
    \sup\limits_h\frac{\omega_1(f,|h|_1)}{|h|_1}=\sup\limits_x\|f'(x)\|_{E_1\rightarrow
    E_2}.
\end{equation*}

При любом $r\in \mathbb{N}$ $\ \omega_r$ убывает вместе с $|h|_1$ и
$$\omega_r(f;|h|_1)\leq2\omega_{r-1}(f;|h|_1)\leq...\leq2^{r-1}\omega_1(f;|h|_1)\leq2^r\|f\|_\infty\equiv2^r\omega_0(f;|h|_1).$$

Исходя из тождества ($n\in \mathbb{N}$)
\begin{equation*}
    \Delta_{n
    \delta}^rf(x)=\sum\limits_{\nu_1=0}^{n-1}...\sum\limits_{\nu_r=0}^{n-1}\Delta_\delta^rf\big(x+(\nu_1+...+\nu_r)\delta\big)\tag{*}
\end{equation*}
(при $r=1$ оно очевидно, а далее доказывается по индукции), получаем

\textbf{4.1.} \emph{$\omega_r(f;n|h|_1)\leq n^r\omega_r(f;|h|_1)$,
$\omega_r(f;\lambda|h|_1)\leq(\lambda+1)^r\omega_r(f;|h|_1)\ $ $(\lambda>0)$.}

\emph{Кроме того, при $|u|_1\leq|v|_1$
\begin{equation*}
    \frac{\omega_r(f;|v|_1)}{|v|_1^r}\leq2^r\frac{\omega_r(f;|u|_1)}{|u|_1^r}
\end{equation*}
(детали см., напр., [\ref{Trigub_Belinsky}--\ref{Dzyadyk_Shevchuk}]).}

\textbf{4.2.} \emph{
\begin{equation*}
    \lim\limits_{h\rightarrow0}\frac{\omega_r(f;|h|_1)}{|h|_1^r}=\sup\limits_h\frac{\omega_r(f;|h|_1)}{|h|_1^r}.
\end{equation*}}

$\rhd$ Для доказательства обозначим верхнюю грань через $M_r\in(0,+\infty]$ и возьмем
произвольное число $M\in(0,M_r)$.

Тогда существует $\delta\neq0$ такое, что $\omega_r(f;|\delta|_1)>M |\delta|_1^r$ и при
$0<|h|_1\leq|\delta|_1$ и $n=\Big[\frac{|\delta|_1}{|u|_1}\Big]$ (целая часть),
учитывая, что $n|h|_1\leq|\delta|_1<(n+1)|h|_1$ и $h\rightarrow0$ при
$n\rightarrow\infty$ получаем из 4.1
\begin{equation*}
\omega_r(f;|\delta|_1)\leq\Big(\frac{|\delta|_1}{|h|_1}+1\Big)^r\omega_r(f;|h|_1)\leq(n+2)^r\omega_r(f;|h|_1).
\end{equation*}

Поэтому
\begin{equation*}
M<\frac{\omega_r(f;|\delta|_1)}{|\delta|_1^r}\leq\frac{\omega_r(f;|\delta|_1)}{n^r|h|_1^r}\leq\Big(1+\frac{2}{n}\Big)^r\frac{\omega_r(f;|h|_1)}{|h|_1^r}\leq\Big(1+\frac{2}{n}\Big)^rM_r.
\end{equation*}
и
\begin{equation*}
    \lim\limits_{h\rightarrow0}\frac{\omega_r(f;|h|_1)}{|h|_1^r}=M_r.
\end{equation*}
$\qquad\qquad\qquad\qquad\qquad\qquad\qquad\qquad\qquad\qquad\qquad\qquad\qquad\qquad\qquad\qquad\qquad\qquad\qquad\blacktriangleleft$

\textbf{4.3.} \emph{При любом целом $k\geq0$
\begin{equation*}
\omega_r(f;|h|_1)\leq\frac{r}{2}\sum\limits_{\nu=0}^k\frac{1}{2^{\nu
r}}\omega_{r+1}(f;2^\nu|h|_1)+\frac{1}{2^{(k+1)r}}\omega_{r}(f;2^{k-1}|h|_1).
\end{equation*}}

$\rhd$ Для доказательства применяем тождество (*) при $n=2$.
\begin{equation*}
     \begin{split}
         & 2^r\Delta_h^rf(x)-\Delta_{2h}^rf(x)=\sum\limits_{\nu=0}^r\binom{r}{\nu}\Big[\Delta_h^rf(x)-\Delta_h^rf(x+\nu h)\Big]= \\
         &=\sum\limits_{\nu=1}^r\binom{r}{\nu}\sum\limits_{m=0}^{\nu-1}\Delta_h^{r+1}f(x+mh)=\sum\limits_{m=0}^{r-1}\Delta_h^{r+1}f(x+mh)\sum\limits_{\nu=m+1}^r\binom{r}{\nu}
     \end{split}
\end{equation*}
(изменен порядок суммирования).

Получаем
\begin{equation*}
    \big|2^r\Delta_h^rf(x)\big|_2\leq\omega_r(f;2|h|_1)+\omega_{r+1}(f;|h|_1)\sum\limits_{m=0}^{r-1}\sum\limits_{\nu=m+1}^r\binom{r}{\nu}.
\end{equation*}

Но двойная сумма равна $r\cdot2^{r-1}$. Поэтому
\begin{equation*}
    \omega_r(f;|h|_1)\leq\frac{r}{2}\omega_{r+1}(f;|h|_1)+\frac{1}{2^r}\omega_r(f;2|h|_1).
\end{equation*}

Заменяя в этом неравенстве $|h|_1$ на $2|h|_1$ и применяя его $k$ раз, получим искомое
соотношение.

Перейдем от сумм к интегралам.

Если $\varphi$ и $\psi$ --- положительные функции на $\mathbb{R}_+=[0,+\infty)$,
$\varphi$ убывает, а $\psi$ возрастает и при некотором $c>0$ $\ \psi(x+1)\leq c\psi(x)$
для $x\in\mathbb{R}_+$, то
\begin{equation*}
    \sum_{\nu=1}^k\varphi(\nu)\psi(\nu)\leq c\int\limits_{0}^k\varphi(x)\psi(x)dx.
\end{equation*}

В нашем случае, когда $\psi(x)=\omega_{r+1}(f;2^x|h|_1)$, $c=2^{r+1}$ (см. 4.1 при
$n=2$). Поэтому
\begin{equation*}
    \omega_r(f;|h|_1)\leq\frac{r}{2}\omega_{r+1}(f;|h|_1)+\frac{1}{2^{(k+1)r}}\omega_r(f;2^{k+1}|h|_1)+2^{r+1}\int\limits_0^k\frac{\omega_{r+1}(f;2^x|h|_1)}{2^{rx}}dx
\end{equation*}
и после замены в интеграле $2^x|h|_1=v$
\begin{equation*}
    \omega_r(f;|h|_1)\leq\frac{r}{2}\omega_{r+1}(f;|h|_1)+\frac{1}{2^{(k+1)r}}\omega_r(f;2^{k+1}|h|_1)+\frac{2^{r+1}}{\ln2}|h|_1^r\int\limits_{|h|_1}^{2^k|h|_1}\frac{\omega_{r+1}(f;v)}{v^{r+1}}dv.
\end{equation*}

При $|h|_1\leq1$, выбирая $k$ так, чтобы $1\leq2^k|h|_1<2$, напр., получим
\begin{equation*}
\begin{split}
    \omega_r(f;|h|_1)&\leq\frac{1}{2^r}\omega_{r}(f;4)|h|_1^r+\frac{r}{2}\omega_{r+1}(f;|h|_1)+\frac{2^{r+1}}{\ln2}|h|_1^r\int\limits_{|h|_1}^{2}\frac{\omega_{r+1}(f;v)}{v^{r+1}}dv\leq\\
    &\leq\frac{1}{2^r}\omega_{r}(f;4)|h|_1^r+c(r)|h|_1^r\int\limits_{|h|_1}^{2}\frac{\omega_{r+1}(f;u)}{u^{r+1}}du.
\end{split}
\end{equation*}
(применено неравенство из 4.1 при $|u|_1=|h|_1$ и $|v|_1\in\big[|h|_1,2\big]$).

Это обобщенное неравенство
Маршо.\qquad\qquad\qquad\qquad\qquad\qquad\qquad\qquad\qquad\qquad$\blacktriangleleft$

\textbf{4.4.} \emph{Если $\omega_{r+1}(f;1)=0$, то при любом $h\in E_1\ $
$\omega_r(f;|h|_1)\equiv\omega_r(f;1)|h|_1^r$.}

$\rhd$ В силу монотонности $\omega_{r+1}(f;|h|_1)=0$ при $|h|_1\leq1$, а значит, и при
всех $h$ (см. 4.1).

В силу 4.3 и 4.1
\begin{equation*}
    \omega_r(f;|h|_1)\leq\frac{1}{2^{(k+1)r}}\omega_r(f;2^{k+1}|h|_1)\leq\omega_r(f;|h|_1).
\end{equation*}

Поэтому при любом $k$ и $h\neq0$
\begin{equation*}
    \frac{\omega_r(f;|h|_1)}{|h|_1^r}=\frac{\omega_r(f;2^{k+1}|h|_1)}{2^{(k+1)r}|h|_1^r}\quad
    \text{или}\quad
    \frac{\omega_r\Big(f;\frac{|h|_1}{2^{k+1}}\Big)}{\Big(\frac{|h|_1}{2^{k+1}}\Big)^r}=\frac{\omega_r(f;|h|_1)}{|h|_1^r}.
\end{equation*}

Осталось перейти к пределу при $k\rightarrow\infty$ (см. 4.2) и положить $|h|_1=1$.
$\qquad\qquad\blacktriangleleft$

\textbf{4.5.} \emph{Для любых $r$ и $m\in\mathbb{N}$ существует $c(r,m)$: при
$|h|_1\leq1$}

\emph{а) }$\omega_{r+m}^r(f;1)\omega^{r+m}_r(f;|h|_1)\leq
c(r,m)\omega^{r+m}_r(f;1)\omega_{r+m}^r(f;|h|_1)$;

\emph{б)} $\omega_{r+m}(f;1)\omega_r(f;|h|_1^{r+m})\leq
c(r,m)\omega_r(f;1)\omega_{r+m}(f;|h|_1^r)$;

\emph{в) }$\omega_{r+m}(f;1)\omega_r(f;|h|_1)\omega_m(f;|h|_1)\leq
c(r,m)\omega_r(f;1)\omega_m(f;1)\omega_{r+m}(f;|h|_1)$.

Доказательство не отличается от приведенного в [\ref{De_Vore}], \textbf{4.6.5}.

\textbf{4.6.} \emph{Если $E_2$ --- банахова алгебра
$\big(|fg|_2\leq|f|_2\cdot|g|_2\big)$, то при $r=1$ и $|h|_1\leq1$
\begin{equation*}
    \omega(fg;|h|_1)\leq\|f\|_\infty\omega(g;|h|_1)+\|g\|_\infty\omega(f;|h|_1),
\end{equation*}
а при $r\geq2$
\begin{equation*}
    \begin{split}
         \omega_r(fg;|h|_1)&\leq c(r)\Big[\|g\|_\infty\cdot\omega_r(f;|h|_1)+\|f\|_\infty\cdot\omega_r(g;|h|_1)+ \\
        &+\|f\|_\infty\cdot\|g\|_\infty\cdot|h|_1^r+\|f\|_\infty\cdot\|g\|_\infty\Big(\frac{\omega_r(f;|h|_1)}{\omega_r(f;1)}+\frac{\omega_r(g;|h|_1)}{\omega_r(g;1)}\Big)\Big],
    \end{split}
\end{equation*}
где по предположению $\frac{0}{0}=0$.}

В [\ref{De_Vore}] такое неравенство доказано для случая $f\in L_p$, $p\in[1,+\infty]$, и
$g\in L_\infty$ (с использованием производных).

$\rhd$ Исходим из того, что при некоторых положительных числовых коэффициентах
$\{a_{k,r}\}$, $\{b_{k,r}\}$ и $\{c_{k,r}\}$ (при $r=1$ сумма отсутствует)
\begin{equation*}
    \Delta_h^r(fg;x)=g(x)\Delta_h^rf(x)+f(x+rh)\Delta_h^rg(x)+\sum\limits_{k=1}^{r-1}a_{k,r}\Delta_{b_{k,r}h}^kf(x)\Delta_{c_{k,r}h}^{r-k}g(x).
\end{equation*}

Следовательно,
\begin{equation*}
    \omega_r(fg;|h|_1)\leq\|g\|_\infty\omega_r(f;|h|_1)+\|h\|_\infty\omega_r(g;|h|_1)+c_1(r)\sum\limits_{k=1}^{r-1}\omega_k(f;|h|_1)\omega_{r-k}(g;|h|_1).
\end{equation*}

При $r\geq2$ предположим сначала, что $\omega_r(f;1)\cdot\omega_r(g;1)\neq0$.

Применяем 4.5 а):
\begin{equation*}
\sum\limits_{k=1}^{r-1}\omega_k(f;|h|_1)\omega_{r-k}(f;|h|_1)\leq
c_2(r)\sum\limits_{k=1}^{r-1}\frac{\omega_k(f;1)}{\omega_r^{\frac{k}{r}}(f;1)}\omega_r^{\frac{k}{r}}(f;|h|_1)\cdot\frac{\omega_{r-k}(g;1)}{\omega_{r-k}^{1-\frac{k}{r}}(g;1)}\omega_r^{1-\frac{k}{r}}(g;|h|_1).
\end{equation*}

Учтем теперь, что $\omega_k(f;1)\leq2^k\|f\|_\infty$,
$\omega_{r-k}(g;1)\leq2^{r-k}\|g\|_\infty$, а при $\delta\in[0,1]$, $x_1>0$ и $x_2>0$
\begin{equation*}
    x_1^\delta x_2^{1-\delta}\leq\max\{x_1,x_2\}<x_1+x_2.
\end{equation*}

Искомое неравенство доказано.

Пусть теперь $\omega_r(f;1)=0$, а $\omega_r(g;1)\neq0$.

Тогда в силу 4.4 $\omega_{r-1}(f;|h|_1)\equiv\omega_{r-1}(f;1)|h|_1^{r-1}$.

При $r\geq3$ применяем неравенство Маршо 4.3 в интегральной форме с заменой $r$ на
$r-2$:
\begin{equation*}
\omega_{r-2}(f;|h|_1)\leq2^r\omega_{r-2}(f;1)|h|_1^{r-2}+c_3(r)\omega_{r-1}(f;1)|h|_1^{r-2}\leq
c_4(r)\omega_{r-2}(f;1)|h|_1^{r-2}.
\end{equation*}

Следовательно, при любом $k\in[1,r-1]$ $\ \omega_{r-k}(f;|h|_1)\leq
c_5(r)\omega_{r-k}(f;1)|h|_1^{r-k}$.

Поэтому
\begin{equation*}
    \sum\limits_{k=1}^{r-1}\omega_k(f;|h|_1)\omega_{r-k}(g;|h|_1)=\sum\limits_{k=1}^{r-1}\omega_{r-k}(f;|h|_1)\omega_k(g;|h|_1)\leq
    c_6(r)\sum\limits_{k=1}^{r-1}\omega_{r-k}(f;1)|h|_1^{r-k}\omega_k(g;|h|_1).
\end{equation*}

В силу 4.1 ($|v|_1=1$)
\begin{equation*}
    |h|_1^{r-k}\leq\frac{2^{r-k}\omega_{r-k}(g;|h|_1)}{\omega_{r-k}(g;1)},
\end{equation*}
а в силу 4.5 в)
\begin{equation*}
\frac{\omega_{r-k}(g;|h|_1)\omega_{k}(g;|h|_1)}{\omega_{r-k}(g;1)}\leq
c_7(r)\frac{\omega_r(g;|h|_1)\omega_{k}(g;1)}{\omega_r(g;1)}.
\end{equation*}

Поэтому
\begin{equation*}
    \sum\limits_{k=1}^{r-1}\omega_k(f;|h|_1)\omega_{r-k}(g;|h|_1)\leq c_8(r)\frac{\omega_r(g;|h|_1)}{\omega_r(g;1)}\sum\limits_{k=1}^{r-1}\omega_{r-k}(f;1)\omega_k(g;1)
\end{equation*}
и можно еще учесть, что $\omega_m(f;|h|_1)\leq2^m\|f\|_\infty$.

Осталося случай $\omega_r(f;1)=\omega_r(g;1)=0$.

В силу предыдущего при $k\in[1,r-1]$
\begin{equation*}
\omega_k(f;|h|_1)\leq c_5(r)\omega_k(f;1)|h|_1^k,\qquad \omega_{r-k}(g;|h|_1)\leq
c_5(r)\omega_{r-k}(g;1)|h|_1^{r-k}.
\end{equation*}

Поэтому
\begin{equation*}
\sum\limits_{k=1}^{r-1}\omega_k(f;|h|_1)\omega_{r-k}(g;|h|_1)\leq
(c_5(r))^2\sum\limits_{k=1}^{r-1}\omega_k(f;1)\omega_{r-k}(g;1)|h|_1^r.
\end{equation*}

Свойство 4.6 доказано.
$\qquad\qquad\qquad\qquad\qquad\qquad\qquad\qquad\qquad\qquad\qquad\qquad\qquad\blacktriangleleft$

\textbf{Замечание 3.} Предположим, что в формуле 4.6 при $\omega_r(g;1)=0$ и
$\omega_r(f;1)\neq0$ можно убрать последнее слагаемое.

Учитывая, что в этом случае $\omega_r(f+g;|h|_1)\equiv\omega_r(f;|h|_1)$, а
$g^2=(f+g)g-fg$, получаем
\begin{equation*}
     \begin{split}
         \omega_r(g^2;|h|_1)&\leq \omega_r((f+g)g;|h|_1)+\omega_r(fg;|h|_1)\leq\\
         &\leq c(r)\Big[\|g\|_\infty\cdot\omega_r(f;|h|_1)+\|f+g\|_\infty\cdot\|g\|_\infty\cdot|h|_1^r+\omega_r(fg;|h|_1)\Big].
     \end{split}
\end{equation*}

Заменяя $f$ на $\varepsilon f$ и переходя к пределу при $\varepsilon\rightarrow0$, имеем
\begin{equation*}
    \omega_r(g^2;|h|_1)\leq c(r)\|g\|_\infty^2\cdot|h|_1^r\qquad(g^2\in W^r).
\end{equation*}

Если считать, что $g$ --- "полином степени не выше $r-1$", то это неравенство типа
Маркова.

Можно изучать и модули гладкости нецелого порядка $r>0$
\begin{equation*}
    \Delta_h^rf(x)=\sum\limits_{\nu=0}^\infty\binom{r}{\nu}(-1)^\nu f(x+\nu h),\qquad
    \binom{r}{\nu}=\frac{r(r-1)...(r-\nu+1)}{\nu!}.
\end{equation*}

Заметим, что
\begin{equation*}
\sum\limits_{\nu=0}^\infty\Big|\binom{r}{\nu}\Big|=\sum\limits_{\nu=0}^{[r]}\binom{r}{\nu}\Big(1+(-1)^{[r]+\nu}\Big).
\end{equation*}

Те же свойства 4.1--4.6 верны для вектор--функций, напр., на шаре $|x|_1\leq1$ (при
естественных ограниченных на величину шага $|h|$).

Функцию, заданную на шаре, можно продолжить с сохранением модуля непрерывности на все
пространство, полагая
\begin{equation*}
    f(x)=f\Big(\frac{x}{|x|_1}(2-|x|_1)\Big)\quad (1<|x|\leq2),\qquad f(x)=f(0)\quad
    (|x|_1>2).
\end{equation*}

По-видимому, можно продолжить функцию с сохранением $\omega_r$ при $r\geq2$ (по порядку)
методом Хестенса (см., напр., [\ref{De_Vore}], [\ref{Trigub_Belinsky}],
\textbf{4.6.12}).

Это были полные модули гладкости.

Частный модуль гладкости (в направлении единичного вектора $e$) определяется разностью
\begin{equation*}
    \underset{e}\Delta_h^rf(x)=\sum\limits_{\nu=0}^r(-1)^\nu\binom{r}{\nu}f(x+\nu
    he)\qquad \big(|h|_1\in(0,1]\big).
\end{equation*}

Смешанная разность ($e_1$ и $e_2$ --- два линейно независимых единичных вектора)
\begin{equation*}
\underset{e_1,e_2}\Delta_{h_1,h_2}^{r_1,r_2}f(x)=\underset{e_2}\Delta_{h_2}^{r_2}\Big(\underset{e_1}\Delta_{h_1}^{r_1}f\Big)(x).
\end{equation*}

Соотношения между полным модулем гладкости, частными и смешанными в случае числовых
функций на $\mathbb{R}^d$ и стандартного базиса см. в [\ref{Timan}]. Эти соотношения
(неравенства) выводятся из некоторых тождеств, доказываемых непосредственно (без
использования преобразования Фурье). Так что те же неравенства справедливы и в случае
банаховых пространств.

\textbf{Один пример.} Полный модуль гладкости второго порядка ведет себя так, как сумма
двух частных модулей того же порядка и смешанного модуля.

В начале этого $\S4$ введена полунорма
\begin{equation*}
    |f|_{W^r}=\sup\limits_h\frac{\omega_r(f;|h|_1)}{|h|_1^r}=\lim\limits_{h\rightarrow0}\frac{\omega_r(f;|h|_1)}{|h|_1}.
\end{equation*}

Положим для ограниченной $f:E_1\rightarrow E_2$ при $\varepsilon>0$
\begin{equation*}
    K_r(\varepsilon,f)=\inf\limits_{g\in W^r}\Big\{|f-g|_2+\varepsilon|g|_{W^r}\Big\}.
\end{equation*}

\textbf{Лемма} \emph{Для того чтобы $K_r(\varepsilon^r,f)\asymp\omega_r(f;\varepsilon)$
необходимо и достаточно, чтобы для любой такой функции $f$ при любом $\varepsilon>0$
существовала $g_\varepsilon(f)\in W^r$ такая, что
\begin{equation*}
    |f-g_\varepsilon(f)|_2\leq c(r)\omega_r(f,\varepsilon),\qquad |g_\varepsilon(f)|_{W^r}\leq
    c(r)\frac{\omega_r(f,\varepsilon)}{\varepsilon^r}.
\end{equation*}
}

И необходимость, и достаточность очевидны.

Отметим еще, что при $E_1=\mathbb{R}$, во всяком случае, можно ввести интеграл как
предел интегральных сумм.

Но тогда можно ввести и функцию типа Стеклова $(h>0)$
\begin{equation*}
    f_{r,h}(x)=\frac{1}{h^r}\int\limits_0^hd\delta_1\int\limits_0^hd\delta_2...\int\limits_0^h\sum\limits_{\nu=1}^r(-1)^{\nu+1}\binom{r}{\nu}f\big(x+\nu\sum\limits_{m=1}^r\delta_m\big)d\delta_r.
\end{equation*}

Ее применение см. в [\ref{Trigub_Belinsky}], \textbf{4.6.7}--\textbf{4.6.10}.

Очевидно, что
\begin{equation*}
    |f-f_{r,h}|_2\leq \omega_r(f,rh),\qquad |f_{r,h}|_{W^r}\leq
    2^{r-1}\frac{\omega_r(f,h)}{h^r},
\end{equation*}
т.к.
\begin{equation*}
    \Delta_h^rf_{r,h}(x)=\frac{1}{h^r}\int\limits_0^hd\delta_1...\int\limits_0^h\sum\limits_{\nu=1}^r(-1)^{\nu+1}\binom{r}{\nu}\Delta_h^rf\big(x+\nu\sum\limits_{m=1}^r\delta_m\big)d\delta_r,
\end{equation*}
\begin{equation*}
    \|f_{r,h}\|_{W^r}\leq\frac{1}{h^r}\sum\limits_{\nu=1}^r\binom{r}{\nu}\omega_r(f;h)=(2^r-1)\frac{\omega_r(f;h)}{h^r}.
\end{equation*}

Да и производную $f^{(r)}$ можно определить как линейный ограниченный оператор
$E_1\rightarrow E_2$ такой, что
\begin{equation*}
    \lim\limits_{h\rightarrow0}\Bigg|\frac{\Delta_h^rf(x)}{|h|_1}-f^{(r)}(x)(h^\circ)^r\Bigg|_2=0\qquad \Big(h^\circ=\frac{h}{|h|_1}\Big).
\end{equation*}

Правда, из существования $f^{(r)}$ при $r\geq2$ не следует существования $f^{(r-1)}$.

\begin{center}
    \textbf{\S5 Некоторые нерешенные вопросы}
\end{center}

1) Указать необходимое и близкое к нему достаточное условие на возрастающую
последовательность $\{\nu_k\}_0^\infty$ натуральных чисел (и, возможно, выпуклую), при
которой средние арифметические частных сумм ряда Фурье
\begin{equation*}
    \frac{1}{n+1}\sum\limits_{k=0}^nS_{\nu_k}(f;x)
\end{equation*}
сходятся при $n\rightarrow\infty$ во всех точках Лебега любой функции $f\in
L_1(\mathbb{T})$ (Zalcwasser, 1936).

Для точек Лебега есть критерий суммируемости, т.е. необходимое и достаточное условие
одновременно ([\ref{Trigub_Belinsky}], \textbf{8.1.3}, см. там же близкие вопросы).

Тот же вопрос для $\gamma$-точек О.~Д.~Габисония (Матем. зам., 1973).

2) Известно, что для $f\in L_p(\mathbb{T}^d)$
\begin{equation*}
    \lim\limits_{n\rightarrow\infty}\Bigg\|f(\cdot)-\sum\limits_{|k|\leq n}\Big(1-\frac{|k|^\alpha}{n^\alpha}\Big)_+^\beta\hat{f}_ke_k\Bigg\|_p=0
\end{equation*}
при $p=1$ и $p=\infty$ ($L_\infty=C$) только в случае $\alpha>0$ и
$\beta>\frac{d-1}{2}$. При $p=2$ ответ очевиден ($\alpha>0$, $\beta\geq0$).

Имеются условия сходимости в зависимости от $\alpha$ и $\beta$ и при некоторых других
$p$ (см. [\ref{Stein}], \textbf{7.8}). Этому вопросу посвящена книга
[\ref{Davis_Chang}]. Но нет полного и точного ответа в общем случае.

3) Как выглядит точный порядок приближения суммами Марцинкевича (средние арифметические
кубических частных сумм Фурье) для периодических функций трех и более переменных (какой
разностный оператор нужно использовать). См. 3.6 выше ($d=2$).

4) Определить специальный модуль непрерывности $\omega^*$ такой, чтобы было для всех
$f\in C(\mathbb{T})$ и $n\in \mathbb{N}$
\begin{equation*}
    \Big\|\frac{1}{n+1}\sum\limits_{k=0}^n\big|f(\cdot)-S_k(f;\cdot)\big|\Big\|_\infty\asymp\omega^*(f;\varepsilon_n)
\end{equation*}
при некоторой нуль-последовательности $\{\varepsilon_n\}$.

5) В этой статье для оценки порядка приближения сверху и снизу сверткой функции с данным
ядром применяется метод мультипликаторов Фурье. Поэтому ответ дается в терминах
коэффициентов ядра, т.е. спектра интегрального оператора. А при каких достаточно общих
условиях на само ядро $K_n$
\begin{equation*}
    \|f-f\ast K_n\|\asymp\omega_r\Big(f;\frac{1}{n}\Big).
\end{equation*}

При $r=1$ и $r=2$ эта задача решается в [\ref{Trigub2009}]. А какой ответ при $r\geq3$?

6) В $\S3$ настоящей статьи изучается приближение операторами эллиптического типа, у
которых символ равен нулю только в нуле.

Следующий (первый) случай изучен Э.~С.~Белинским [\ref{Belinsky1993}] (см. также [\ref{Belinsky1996}]):

Пусть $D_r=\frac{\partial^{r_1+...+r_d}}{\partial x_1^{r_1}...\partial x_d^{r_d}}$ ---
смешанная производная порядка $r=\sum\limits_{j=1}^dr_j$, где
$0<r_1=...=r_\nu<r_{\nu+1}\leq...\leq r_d$. Предположим, что $\int\limits_{-\pi}^\pi
f(x)dx_j=0$ ($1\leq j\leq d$) и для $s\in \mathbb{Z}_+^d$
\begin{equation*}
    \delta_s(f)=\sum\limits_{k\in \rho(s)}\hat{f}_ke_k,\qquad
    \rho(s)=\big\{k\in\mathbb{Z}^d: 2^{s_j-1}\leq|k_j|<2^{s_j},\ 1\leq j\leq d\big\}.
\end{equation*}

Тогда при $p\in(1,+\infty)$ и $\Big[\frac{1}{\varepsilon}\Big]=2^n$ ($J$ --- единичный
оператор)
\begin{equation*}
    K\big(\varepsilon^{r_1},f,L_p,W(D_r)_p\big)\asymp\Big\|f-\sum\limits_{s:(s,r)\leq r,n}\big(J-2^{-r,n}D_r\big)\delta_s(f)\Big\|_p.
\end{equation*}

Каков модуль гладкости, порожденный смешанной производной в $L_p$ ($1\leq p\leq\infty$),
и $K$-функционал в $C$?

7) Для функций на отрезке роль концевых точек ососбая.

Два важных результата о Д-Т модулях см. в замечании 1 ($\S2$).

Приведем еще один (V.~Totik [\ref{Totik}], см. также [\ref{De_Vore}]):

для классических полиномов Бернштейна $B_n$ на $[-1,1]$
\begin{equation*}
    \|f-B_n(f)\|_C\asymp\omega_2^\varphi\Big(f;\frac{1}{\sqrt{n}}\Big).
\end{equation*}

Но существенно ранее были хорошо изучены приближения алгебраическими полиномами (прямые
и обратные теоремы) с учетом положения точки (С.~М.~Никольский, А.~Ф.~Тиман,
В.~К.~Дзядык, G.~Freud, Ю.~А.~Брудный). См., напр., [\ref{De_Vore}],
[\ref{Trigub_Belinsky}], \textbf{4.7}. А для тех же полиномов Бернштейна давно доказано,
что на $[-1,1]$
\begin{equation*}
    |f(x)-B_n(f,x)|\leq
    c\omega_2\Big(f;\Big(\frac{\sqrt{1-x^2}}{n}\Big)^{\frac{1}{2}}\Big)_\infty,
\end{equation*}
а такая же оценка приближения снизу не верна [\ref{Cao_Gonska_Kasco}].

Для функций на отрезке есть ряды Фурье-Якоби, для функций на множествах комплексной
плоскости --- ряды Фабера. Как сформулировать и применить принцип сравнения методов
суммирования таких рядов (с учетом положения точки)?

Возможны ли двусторонние оценки приближения на $[-1,1]$ вида
\begin{equation*}
    \omega_r\Big(f;\frac{\sqrt{1-x^2}}{n}+\frac{1}{n^2}\Big)_\infty?
\end{equation*}

В случае нелинейных операторов приближения возможны:

Пусть $s\in \mathbb{Z}_+$. Тогда для любой функции $f\in C^s[-1,1]$ при $r\in\mathbb{N}$
и $n\geq\max(r+s-1,3)$ существует полином $p_n$ степени не выше $n$ такой, что при
$x\in[-1,1]$
\begin{equation*}
    \delta_n^s(x)\omega_r\big(f^{(s)};\delta_n(x)\big)\leq p_n(x)-f(x)\leq
    c(s,r)\delta_n^s(x)\omega_r\big(f^{(s)};\delta_n(x)\big),
\end{equation*}
где $\delta_n(x)=\frac{\sqrt{1-x^2}}{n}+\frac{1}{n^2}$ (см. [\ref{Trigub2000}]).

\end{document}